%% file: sfPIF4_arXiv.tex
\documentclass[times,preprint,3p]{elsarticle}

\journal{Journal of Computational Physics}

\usepackage{framed,multirow}

\input{header}

\begin{document}
\begin{frontmatter}

\title{A recursive system-free single-step temporal discretization method 
for finite difference methods}

\author[1]{Youngjun {Lee}}
\ead{ylee109@ucsc.edu}
\author[1]{Dongwook {Lee}\corref{cor1}}
\cortext[cor1]{Corresponding author: Tel.: +1-831-502-7708}
\ead{dlee79@ucsc.edu}
\author[2]{Adam {Reyes}}
\ead{adam.reyes@rochester.edu}
\address[1]{Department of Applied Mathematics, The University of California, Santa Cruz, CA, United States}
\address[2]{Department of Physics and Astronomy, The University of Rochester, NY, United States}

\begin{abstract}
Single-stage or single-step high-order temporal discretizations of partial differential equations (PDEs) have shown great promise in delivering high-order accuracy in time with efficient use of computational resources. There has been much success in developing such methods for finite volume method (FVM) discretizations of PDEs. The Picard Integral formulation (PIF) has recently made such single-stage temporal methods accessible for finite difference method (FDM) discretizations.  PIF methods rely on the so-called Lax-Wendroff procedures to tightly couple spatial and temporal derivatives through the governing PDE system to construct high-order Taylor series expansions in time. Going to higher than third order in time requires the calculation of \textit{Jacobian-like} derivative tensor-vector contractions of an increasingly larger degree, greatly adding to the complexity of such schemes. To that end, we present in this paper a method for calculating these tensor contractions through a recursive application of a discrete Jacobian operator that readily and efficiently computes the needed contractions entirely agnostic of the system of partial differential equations (PDEs) being solved.
\end{abstract}

\begin{keyword}
    Picard integration formulation;
    Jacobian-free and Hessian-free;
    Recursive;
    Cauchy-Kowalewski procedure;
    High-order method;
    Finite difference method
\end{keyword}
\end{frontmatter}



\section{Introduction}\label{sec:introduction}

In past decades, high-order discrete methods for hyperbolic conservation laws
have become one of the central themes in computational fluid dynamics (CFD)
due to their potential in achieving highly accurate predictions
in balance with trends in current high-performance computing (HPC) architectures.
Using high-order methods allows better solution accuracy
with a fixed memory profile at the cost of increased floating point operations;
in line with next generation HPC systems that are seeing increased 
computing power with saturation in the memory per compute core. 

Practitioners have focused on designing
high-arithmetic-intensity data reconstruction and interpolation models
in the context of finite volume methods (FVM) and finite difference methods (FDM).
Under the dual computational need for accuracy and stability,
the CFD community has developed high-order reconstruction and interpolation methods
that can produce highly accurate solutions while avoiding unphysical
oscillations near the discontinuities.
Examples include the early success of the piecewise parabolic method (PPM)
by Colella and Woodward~\cite{colella1984piecewise,woodward1984numerical},
which has been still actively adopted as a shock-capturing partial differential equation (PDE)  
solver by many CFD users after about four decades
since its introduction.
In 1987, Harten \textit{et al.}~\cite{harten1987uniformly}
proposed an essentially non-oscillatory (ENO) scheme
that chooses the appropriate stencil adaptively
to prevent spurious oscillations near strong gradients.
Liu \textit{et al.}~\cite{liu1994weighted} improved this idea
by introducing the weighted essentially non-oscillatory (WENO) scheme,
which takes a convex combination of all possible stencils,
reducing the number of logical calculations needed for the original ENO scheme.
The WENO scheme was further improved by Jiang and Shu~\cite{jiang1996efficient}
and became one of the most popular high-order reconstruction and interpolation methods
for solving shock-dominant CFD simulation.
Several modifications of WENO methods have been proposed over decades,
such as WENO-Z~\cite{borges2008improved, castro2011high},
central-WENO~\cite{levy1999central, qiu2002construction},
Hermite-WENO~\cite{qiu2004hermite,qiu2005hermite}
WENO-AO~\cite{balsara2016efficient}, 
and polynomial-free GP-WENO~\cite{reyes2018new,reyes2019variable},
to name a few.

Common to all numerical PDE solvers is the solution lies on the spatio-temporal plane.
Numerical errors arise from both the spatial and temporal axis;
thereby, a high-order scheme requires a meticulously designed
temporal discretization method that ensures the solution's accuracy and stability.
However, efforts to achieve a high-order of accuracy in the temporal axis have seen a renewed effort. 
For decades, the strong stability preserving Runge-Kutta (SSP-RK)
time integrator has been considered as the standard temporal integration strategy
for an extensive range of high-order numerical schemes for PDE solvers.
The key idea of SSP-RK scheme is to maintain the strong stability property (SSP)
at high-order accuracy by sequentially applying convex combinations
of the first-order forward Euler method as a building-block
at each sub-stage of the Runge-Kutta method.
In this manner, the total variation diminishing (TVD) property
is achieved in each sub-stages; therefore, it ensures the whole scheme's stability.

Despite the high portability and fidelity of SSP-RK methods,
the very nature of the SSP-RK method~-- being a multi-stage approach~--
increases computational resources in CFD simulations.
In SSP-RK methods, the data reconstruction/interpolation
and the boundary condition should be applied in each sub-stage,
which increases the computational costs and the footprint of data communications
in the parallel computational architecture.
It makes the simulations using
the adaptive mesh refinement (AMR) method less attractive,
which progressively refines the grid resolutions and increases data communications
around the simulations' interesting features.

To alleviate such issues, many practitioners recently have proposed single-step,
high-order time updating strategies based on the so-called Lax–Wendroff (Taylor) method.
The core design principle of the Lax-Wendroff method is to convert
time derivatives to spatial derivatives through
the Lax-Wendroff or Cauchy-Kowalewski procedure (LW/CK hereafter).
In this way, the spatial and temporal derivatives are coupled through
Jacobians and Hessians; thus, the numerical solutions can be updated
in a single step while maintaining the high-order accuracy.
In 2001, Toro \textit{et al.}~\cite{toro2001towards} extended this idea
by combining it with the generalized Riemann problems (GRPs) and introduced 
the Arbitrary high order DErivative Riemann problem (ADER) method.
Toro and his collaborators applied LW/CK procedures
to get the coefficients of the power series expansion
of the conservative variables and solved GRPs for each high-order term.
ADER methods were further developed in~\cite{toro2001towards,titarev2002ader,titarev2005ader},
and it has grown its popularity over decades,
leading to various further modifications, including
ADER-DG~\cite{fambri2017space,zanotti2016efficient},
ADER-CG~\cite{balsara2009efficient,balsara2013efficient,balsara2017higher},
and ADER-Taylor~\cite{norman2012multi,norman2013algorithmic,norman2014weno},
to name a few.
Comparisons with RK methods have shown that these
single-stage time updates are more efficient in terms of CPU time to solution~\cite{balsara2013efficient,lee2017piecewise,lee2021single}.

The Picard integral formulation (PIF) method,
proposed by Christlieb \textit{et al.}~\cite{christlieb2015picard},
is another Lax-Wendroff type time discretization method
under the finite difference formulation, which evolves the pointwise conservative variables.
Instead of taking temporally pointwise numerical fluxes as in the conventional FDM,
the PIF method takes time-averaged numerical fluxes for
updating the solutions in a single step. Christlieb and his collaborators
demonstrated that LW/CK procedures could successfully be utilized
for obtaining high-order terms of the numerical fluxes
as in many Lax-Wendroff type works of literature,
which results in faster CPU performance than the SSP-RK method with the same order of accuracy.

The primary advantage of LW/CK-based methods is the enhanced performance offered by being a single-step method, harnessing the tight coupling of temporal and spatial derivatives through LW/CK procedures to construct high-order time-series Taylor expansions.
This becomes hugely attractive in massively parallel computing, minimizing the computational frequency of data transfers between processors each time step, which would need to be repeated for each intermediate RK stage. On the other hand, the dependence of the strong coupling on analytic derivatives of the governing PDEs makes the LW/CK approach 
less flexible and less broadly applicable to all systems of PDEs.

To meet this end, in~\cite{lee2021single}, we proposed a novel approach
to bypass the analytic evaluations of Jacobian and Hessian terms,
which are the major implementation hurdles for Lax-Wendroff type methods.
This new approach, which we called the system-free (SF) method,
adopts the Jacobian-free idea commonly used for
iterative methods~\cite{gear1983iterative,brown1990hybrid,knoll2004jacobian,knoll2011application}.
We implemented this new approach to the original third-order PIF method,
called the third-order SF-PIF method (or SF-PIF3).
The new SF-PIF3 method shows the same performance results while maintaining
the same order of accuracy and stability as the original PIF method.
Also, by virtue of the system-independent approach for the \textit{Jacobian-like} tensor contractions,
SF-PIF3 can apply to a different hyperbolic system of equations with ease.

The original PIF method and SF-PIF3 are third-order in time,
coupled with the fifth-order spatial reconstruction method (WENO-JS).
Assuming that the temporal integration is discretized with a \( q \)th order scheme
and the spatial with a \( p \)th order scheme (often with \( q \le p \)),
the overall solution accuracy is determined by the leading errors of
the two discretizations, namely, \( \mathcal{O} \left( \Delta t^{q} , \Delta s^{p} \right) \),
where \( \Delta t \) and \( \Delta s \) represent the temporal and spatial length scales.
Practically speaking, the spatial errors usually dominate the temporal errors; thus,
the high-order methods with \( q \le p \) are justifiable.
Nonetheless, we observe that the temporal errors gradually dominate
the spatial errors as we increase the grid resolution progressively,
which leads us to find a higher than a third-order temporal scheme
for maintaining overall solution accuracy even in the finer grid resolutions.

In Lax-Wendroff type time discretization methods, like PIF and SF-PIF3 methods,
rank-4 \textit{Jacobian-like} tensors are needed for the fourth-order approximation.
Although possible, the need for such analytical handling becomes prohibitive in 
the Lax-Wendroff type schemes when extending their accuracy beyond third-order
due to the drastic growth in complexity.

In this regard, we propose a new fourth-order extension of the SF-PIF3 method.
Since SF-PIF3 bypasses all the \textit{Jacobian-like} calculations,
our SF approach is further rewarded in a fourth-order extension, which demands
a leap in calculation counts for conventional Lax-Wendroff type schemes.
Moreover, we present a new improved version of the previous SF approach~\cite{lee2021single},
which applies \textit{a recursive strategy} to obtain the higher-order derivative tensor-vector contractions,
promising a more compact code structure and faster performance than the original SF method.
With such modification, our fourth-order SF-PIF4 method shows
nearly twice faster performance than the optimal fourth-order five-stage SSP-RK method.

We organize the paper as follows. In \cref{sec:pif} we briefly review
the general discretization strategy of the original PIF method.
We present the fourth-order extension of the original PIF method
in \cref{sec:pif4}, and we apply the original SF approach to the fourth-order
PIF method in \cref{sec:original-sf}.
The improved recursive SF approach will be introduced in \cref{sec:new-sf},
which ensures faster performance and a more straightforward  code structure.
The results of various 2D and 3D benchmark problems are presented in \cref{sec:results}.
We conclude our paper with a summary in \cref{sec:conclusion}.

\section{Picard integral formulation}\label{sec:pif}

We are interested in solving the general conservation system of equations in 3D,
\begin{equation}\label{eq:gov}
    \partial_{t} \bU + \nabla \cdot \mathcal{F}(\bU)
    = \partial_{t} \bU + \partial_{x} \bF (\bU) + \partial_{y} \bG (\bU) + \partial_{z} \bH (\bU) = 0,
\end{equation}
where \( \bU \) is the vector of conservative variables
and \( \bF, \bG \), and \( \bH \) are the flux functions in
\( x \)-, \( y \)- and \( z \)-direction, respectively.
In the Euler equations, the conservative variables and
the flux functions are defined as,
\begin{equation}\label{eq:euler-3d}
    \bU = \begin{bmatrix}
        \rho \\
        \rho u \\
        \rho v \\
        \rho w \\
        E
    \end{bmatrix},\quad
    \bF (\bU) = \begin{bmatrix}
        \rho u \\
        \rho u^{2} + p \\
        \rho u v \\
        \rho u w \\
        u \left( E + p \right)
    \end{bmatrix}, \quad
    \bG (\bU) = \begin{bmatrix}
        \rho v \\
        \rho u v \\
        \rho v^{2} + p \\
        \rho v w \\
        v \left( E + p \right)
    \end{bmatrix}, \quad
    \bH (\bU) = \begin{bmatrix}
        \rho w \\
        \rho u w \\
        \rho v w \\
        \rho w^{2} + p \\
        w \left( E + p \right)
    \end{bmatrix}.
\end{equation}

Applying the Picard integral formulation (PIF)~\cite{seal2014high},
we take a time average of~\cref{eq:gov} within a single time step
\( \Delta t \) over an interval
\( [ t^{n}, t^{n} + \Delta t ] = [ t^{n}, t^{n + 1} ] \),
\begin{equation}\label{eq:semi-discrete}
    \bU^{n + 1} = \bU^{n} - \Delta t \left( \partial_{x} \bF^{avg} + \partial_{y} \bG^{avg} + \partial_{z} \bH^{avg} \right),
\end{equation}
where \( \bF^{avg}, \bG^{avg} \), and \( \bH^{avg} \) are the time-averaged fluxes
in each direction, 
\begin{equation}\label{eq:avg-flx}
        \bF^{avg} (\bx) \equiv \frac{1}{\dt} \int^{t^{n + 1}}_{t^{n}} \bF(\bU(\bx, t)) \mathop{dt}, \quad
        \bG^{avg} (\bx) \equiv \frac{1}{\dt} \int^{t^{n + 1}}_{t^{n}} \bG(\bU (\bx, t)) \mathop{dt}, \quad
        \bH^{avg} (\bx) \equiv \frac{1}{\dt} \int^{t^{n + 1}}_{t^{n}} \bH(\bU (\bx, t)) \mathop{dt},
\end{equation}
for $\bx = (x, y, z) \in \mathbb{R}^{3}$.

We wish to express the spatial derivatives of the time-averaged fluxes
in~\cref{eq:semi-discrete}
using highly approximated numerical fluxes
\( \hat{\bff}, \hat{\bg} \), and \( \hat{\bh} \)
at cell interfaces. Taking \( x \)-directional
flux \( \bF \), for example, we aim to express
\begin{equation}\label{eq:flx-deriv}
    \partial_{x} \bF^{avg} \bigg|_{\bx = \bx_{ijk}} =
    \frac{1}{\dx} \left( \hat{\bff}_{i + \half, j, k} -
                                \hat{\bff}_{i - \half, j, k} \right) +
                                \mathcal{O} (\dx^{p} + \dt^{q}), \quad \quad
    \bx_{ijk} = (x_i, y_j, z_k).
\end{equation}
The \( y \)- and \( z \)- directional derivatives are approximated
in a similar fashion.

Finding approximated solutions for the numerical fluxes in the PIF method
is nearly identical to the conventional finite difference method (FDM).
Following the standard convention in FDM, we treat pointwise \( x \)-directional
flux \( \bF(x, y_{j}, z_{k}) \) as a 1D cell average of an auxiliary function
\( \hat{\bF} \) in 1D,
\begin{equation}\label{eq:aux-func}
    \bF(x, y_{j}, z_{k}) = \frac{1}{\Delta x} \int^{x + \frac{\Delta x}{2}}_{x + \frac{\Delta x}{2}}
    \hat{\bF} (\xi, y_{j}, z_{k}) \mathop{d\xi}.
\end{equation}
Then, by taking an analytical derivative of~\cref{eq:aux-func}
with respect to \( x \), we have
\begin{equation}\label{eq:aux-deriv}
    \partial_{x} \bF \bigg|_{\bx = \bx_{ijk}} = \frac{1}{\Delta x} \left( \hat{\bF} (x_{i + \half}, y_{j}, z_{k})
        - \hat{\bF} (x_{i - \half}, y_{j}, z_{k}) \right).
\end{equation}
By comparing~\cref{eq:flx-deriv} and~\cref{eq:aux-deriv},
we can identify the auxiliary function \( \hat{\bF} \) as
a numerical flux \( \hat{\bff} \). We can follow
similar procedures to identify the rest of numerical fluxes \( \hat{\bg} \) and \( \hat{\bh} \).

Therefore, determining \textit{spatially} approximated solutions
for the numerical fluxes is the inverse problem of~\cref{eq:aux-func}, viz. finding
pointwise values at cell interfaces, given
the volume-averaged quantities at cell centers, which
is exactly the same way as the high-order reconstruction procedure
used in 1D finite volume method (FVM).
Namely, we can use the conventional high-order reconstruction
procedures used in FVM for approximating numerical fluxes
\( \hat{\bff}, \hat{\bg} \) and \( \hat{\bh} \),
at the desired \textit{spatially} \( p \)th-order accuracy
by using the cell-centered pointwise analytic fluxes,
\( \bF_{ijk}, \bG_{ijk} \), and \( \bH_{ijk} \).
One strategic difference in the PIF method is that we use the \textit{time-averaged fluxes},
\( \bF^{avg}, \bG^{avg} \), and \( \bH^{avg} \),
as the inputs of the high-order reconstruction method
instead of the pointwise flux values in the conventional FDM
in order to assure the anticipated
\( q \)th-order temporal accuracy
for the numerical fluxes \( \hat{\bff}, \hat{\bg} \), and \( \hat{\bh} \)
at cell interfaces.

The time-averaged fluxes
are obtained through the Taylor expansion of the pointwise
flux around \( t^{n} \).
In the \( q \)th-order PIF method,
the time-averaged \( x \)-directional flux \( \bF^{avg} \)
is approximated as,
\begin{equation}\label{eq:flx-taylor}
    \begin{split}
        \bF^{avg} (\bx)
        &= \frac{1}{\dt} \int^{t^{n + 1}}_{t^n} \bF(\bx, t) \mathop{dt}\\
        &= \bF (\bx, t^{n})
            + \left. \frac{\Delta t}{2!} \partial_{t}^{(1)} \bF (\bx, t) \right|_{t = t^{n}}
            + \left. \frac{\Delta t^{2}}{3!} \partial_{t}^{(2)} \bF (\bx, t) \right|_{t = t^{n}}
            + \left. \frac{\Delta t^{3}}{4!} \partial_{t}^{(3)} \bF (\bx, t) \right|_{t = t^{n}}
            + \cdots \\
        &= \sum\limits_{i=0}^{q-1}
            \left.\frac{\dt^{i}}{(i+1)!} \partial_{t}^{(i)} \bF (\bx, t) \right|_{t = t^{n}} + \mathcal{O}(\Delta t^{q}) \\
        &= \bF^{appx,q} (\bx, t^{n}) + \mathcal{O}(\Delta t^{q}).
    \end{split}
\end{equation}
We will use the \textit{temporally} \( q \)th-order approximated fluxes \( \bF^{appx, q} \)
as the inputs of the \( p\)th-order reconstruction scheme \( \mathcal{R}(\cdot) \)
that is combined with a characteristic flux splitting method \( \mathcal{FS}(\cdot) \) 
to apply the \( p \)th-order \textit{spatial} approximation to
the numerical flux \( \hat{\bff} \) at cell interfaces,
\begin{equation}\label{eq:num-flx}
    \hat{\bff}_{i + \half, j, k} =
        \mathcal{R}\left(
            \mathcal{FS}\left(\bF^{appx, q}_{i-r, j}, \dots,
                \bF^{appx, q}_{i+r+1, j}
        \right)
    \right)
    + \mathcal{O}(\dx^{p}),
\end{equation}
where \( r \) represents the stencil radius
required for the \( p \)th-order reconstruction method, \( \mathcal{R}(\cdot) \).
This study uses the conventional fifth-order WENO-JS
method~\cite{jiang1996efficient} and
the global Lax–Friedrichs flux splitting scheme taking the maximum wave speed
over the entire domain and all characteristic fields. The choice of this global Lax-Friedrichs
scheme is particularly more diffusive than other possible forms of splitting,
although we have found that the added numerical dissipation becomes
less significant for designing our fourth-order SF-PIF scheme without sacrificing accuracy.

Recall that the numerical flux in~\cref{eq:num-flx}
is a high-order approximated solution \textit{both in time and space}
since we used temporally approximated inputs, \( \bF^{appx, q} \),
instead of the temporally pointwise flux values.
Therefore, the \( x \)-directional numerical flux
derived from~\cref{eq:flx-taylor} and~\cref{eq:num-flx}
allows us to update the solution in a single step
while maintaining high order accuracy both in space and time.
The \( y \)- and \( z \)-directional numerical fluxes,
\( \hat{\bg} \) and \( \hat{\bh} \), are obtained in similar fashions.

Finally, the fully discretized form of the governing equation is given as,
\begin{equation}\label{eq:gov-dis}
    \bU^{n + 1}_{i,j,k} = \bU^{n}_{i, j, k} -
    \frac{\dt}{\dx} \left( \hat{\bff}_{i + \half, j, k} -
                                 \hat{\bff}_{i - \half, j, k} \right) -
    \frac{\dt}{\dy} \left( \hat{\bg}_{i, j + \half, k} -
                                  \hat{\bg}_{i, j - \half, k} \right) -
    \frac{\dt}{\dz} \left( \hat{\bh}_{i, j, k + \half} -
                                  \hat{\bh}_{i, j, k - \half} \right).
\end{equation}

\section{The fourth-order extension of the PIF method}\label{sec:pif4}

The primary goal in this section is to extend the PIF method~\cite{seal2014high}
to fourth-order in time. We consider a fourth-order approximated
time-averaged flux in \( x \)-direction, \( \bF^{appx, 4} \),
from the Taylor expansion of the pointwise flux around \( t^{n} \).
As expressed in~\cref{eq:flx-taylor} we have,
\begin{equation}\label{eq:flx-taylor4}
    \bF^{appx, 4} (\bx)
    = \bF (\bx, t^{n})
        + \left. \frac{\Delta t}{2!} \partial_{t}^{(1)} \bF (\bx, t) \right|_{t = t^{n}}
        + \left. \frac{\Delta t^{2}}{3!} \partial_{t}^{(2)} \bF (\bx, t) \right|_{t = t^{n}}
        + \left. \frac{\Delta t^{3}}{4!} \partial_{t}^{(3)} \bF (\bx, t) \right|_{t = t^{n}}.
\end{equation}
The other \( y \)- and \( z \)-directional approximated fluxes,
\( \bG^{appx, 4} \) and \( \bH^{appx, 4} \),
are defined in similarly.
Our main interest is transforming all the
time derivatives in~\cref{eq:flx-taylor4} to the corresponding spatial derivatives;
thereby we could express~\cref{eq:flx-taylor4} in a fully explicit form.

For simplicity, we begin to adopt a compact subscript notation of partial derivatives
and omit
the temporal expression of \( t = t^{n} \). Thus, we rewrite
\cref{eq:gov} in a compact form as,
\begin{equation}\label{eq:gov-compact}
    \bU_{t} + \nabla \cdot \mathcal{F}(\bU) = \bU_{t} + \bF_{x} + \bG_{y} + \bH_{z} = 0.
\end{equation}

In~\cite{lee2021single}, we introduced the so-called
flux equation~--~an evolution equation of fluxes -- 
by applying a chain rule to~\cref{eq:gov-compact}. 
For example, the flux equation of the $x$-flux is
\begin{equation}\label{eq:flx-eq}
    \bF_{t} + \bF_{\bU} \cdot \Div = 0, \quad \text{where}\;\; \Div = \bF_{x} + \bG_{y} + \bH_{z}.
\end{equation}
The above flux equation allows us to convert
time derivatives of the fluxes to the spatial derivatives
via the LW/CK procedure.
For instance, the first-order time derivative of \( \bF \) is
easily converted to the spatial derivatives as,
\begin{equation}\label{eq:Ft}
    \bF_{t} = - \bF_{\bU} \cdot \Div.
\end{equation}
The higher-order time derivatives could be achieved
by taking partial derivatives to~\cref{eq:Ft} recursively.
As an example, the second-order term is written as
\begin{equation}\label{eq:Ftt}
    \bF_{tt} = \bF_{\bU \bU} \cdot \Div \cdot \Div - \bF_{\bU} \cdot \Div_{t},
\end{equation}
where
\begin{equation}\label{eq:divt}
        \Div_{t} = 
         -\bF_{\bU \bU} \cdot \bU_{x} \cdot \Div - \bF_{\bU} \cdot \Div_{x} 
         -\bG_{\bU \bU} \cdot \bU_{y} \cdot \Div - \bG_{\bU} \cdot \Div_{y}
 	 -\bH_{\bU \bU} \cdot \bU_{z} \cdot \Div - \bH_{\bU} \cdot \Div_{z}.
\end{equation}

Following the same procedure, we are able to obtain
an explicit form of the third-order time derivative
of the flux as,
\begin{equation}\label{eq:Fttt}
    \bF_{ttt} = -\bF_{\bU \bU \bU} \cdot \Div \cdot \Div \cdot \Div
    + 3 \bF_{\bU \bU} \cdot \Div \cdot \Div_{t}
    - \bF_{\bU} \cdot \Div_{tt},
\end{equation}
where
\begin{equation}\label{eq:divtt}
    \begin{split}
        \Div_{tt} = \quad & \bF_{\bU \bU \bU} \cdot \Div \cdot \bU_{x} \cdot \Div + 2 \bF_{\bU \bU} \cdot \Div \cdot \Div_{x} - \bF_{\bU \bU} \cdot \bU_{x} \cdot \Div_{t} - \bF_{\bU} \cdot \Div_{tx} \\
            + & \bG_{\bU \bU \bU} \cdot \Div \cdot \bU_{y} \cdot \Div + 2 \bG_{\bU \bU} \cdot \Div \cdot \Div_{y} - \bG_{\bU \bU} \cdot \bU_{y} \cdot \Div_{t} - \bG_{\bU} \cdot \Div_{ty} \\
            + & \bH_{\bU \bU \bU} \cdot \Div \cdot \bU_{z} \cdot \Div + 2 \bH_{\bU \bU} \cdot \Div \cdot \Div_{z} - \bH_{\bU \bU} \cdot \bU_{z} \cdot \Div_{t} - \bH_{\bU} \cdot \Div_{tz},
    \end{split}
\end{equation}
and
\begin{equation}\label{eq:divtx}
    \begin{split}
        \Div_{tx} = \quad & \bF_{\bU \bU \bU} \cdot \bU_{x} \cdot \Div \cdot \bU_{x} - \bF_{\bU \bU} \cdot \bU_{xx} \cdot \Div
        -2\bF_{\bU \bU} \cdot \Div_{x} \cdot \bU_{x} - \bF_{\bU} \cdot \Div_{xx} \\
            - & \bG_{\bU \bU \bU} \cdot \bU_{x} \cdot \Div \cdot \bU_{y} - \bG_{\bU \bU} \cdot \bU_{xy} \cdot \Div - \bG_{\bU \bU} \cdot \Div_{x} \cdot \bU_{y}
        -\bG_{\bU \bU} \cdot \bU_{x} \cdot \Div_{y} - \bG_{\bU} \cdot \Div_{xy} \\
            - & \bH_{\bU \bU \bU} \cdot \bU_{x} \cdot \Div \cdot \bU_{z} - \bH_{\bU \bU} \cdot \bU_{xz} \cdot \Div - \bH_{\bU \bU} \cdot \Div_{x} \cdot \bU_{z}
        -\bH_{\bU \bU} \cdot \bU_{x} \cdot \Div_{z} - \bH_{\bU} \cdot \Div_{xz},
    \end{split}
\end{equation}
and similarly for $\Div_{ty}$ and $\Div_{tz}$.
Collecting~\crefrange{eq:Ft}{eq:divtx}, we can express the fourth-order
approximation of the time-averaged flux \( \bF^{appx,4} \) explicitly,
as the spatial derivatives are readily approximated through the conventional
central differencing schemes.

We should note that we adopt the conventional five-points, fourth-order in space
central differencing schemes to evaluate the spatial derivative terms,
following the original PIF method in~\cite{seal2014high, christlieb2015picard}.
Moreover, for reducing the code complexity and improving the code performance,
we reuse the divergence of the flux, \( \Div \), for calculating
high order spatial derivatives, e.g., \( \Div_{x}, \Div_{xx} \), and \( \Div_{xy} \).
This approach requires an additional guard cell layer
(resulting in two more guard cells for the five-points derivatives). However,
the overall code performance is better than evaluating
high-order derivatives in each direction.
Also, we have observed that it does not affect the accuracy and the stability
of the PIF scheme.

\section{The original non-recursive system-free (SF) approach}\label{sec:original-sf}

As firstly proposed in~\cite{lee2021single},
the system-free (SF) method provides a capability to bypass
all the analytical derivations of
\textit{Jacobian-like} terms (\( \bF_{\bU}, \bF_{\bU\bU}, \cdots \))
in the PIF method by considering a central differencing
with a small perturbation \( \varepsilon \)
for the input space of the flux function. For example,
a dot product between the Jacobian matrix \( \bF_{\bU} \)
and an arbitrary vector \( \bV \) can be approximated as,
\begin{equation}\label{eq:orig-sf-jac}
    \bF_{\bU} \cdot \bV = \frac{1}{2\varepsilon}
    \bigg[ \bF(\bU + \varepsilon \bV) -\bF(\bU - \varepsilon \bV) \bigg]
    + \mathcal{O} ( \varepsilon^{2} ).
\end{equation}
We can extend the above idea to the Hessian tensor contraction with the two same vectors \( \bV \),
\begin{equation}\label{eq:orig-sf-hes}
    \bF_{\bU \bU} \cdot \bV \cdot \bV = \frac{1}{\varepsilon^{2}}
    \bigg[ \bF(\bU + \varepsilon \bV) -2\bF(\bU) -\bF(\bU - \varepsilon \bV) \bigg]
    + \mathcal{O} ( \varepsilon^{2} ).
\end{equation}
The contraction with two different vectors 
of \( \bV \) and \( \bW \)
is achieved by a linear combination of~\cref{eq:orig-sf-hes}
by following the simple vector calculus,
\begin{equation}\label{eq:orig-sf-hes-vec}
    \bF_{\bU \bU} \cdot \bV \cdot \bW = \frac{1}{2}
    \bigg[ \bF_{\bU \bU} \cdot \left( \bV + \bW \right) \cdot \left( \bV + \bW \right) -
          \left( \bF_{\bU \bU} \cdot \bV \cdot \bV + \bF_{\bU \bU} \cdot \bW \cdot \bW \right) \bigg].
\end{equation}

Theoretically speaking, the original system-free procedure in the above can be
applied to any arbitrary order of derivatives of the flux function \( \bF \)
with respect to the conservative variable \( \bU \). For instance,
the fourth-order extension of the PIF method
requires the third-order derivative of \( \bF \), i.e., \( \bF_{\bU \bU \bU} \).
Following the same mathematical basis of~\cref{eq:orig-sf-jac}
and~\cref{eq:orig-sf-hes}, we obtain
\begin{equation}\label{eq:orig-sf-don}
    \bF_{\bU \bU \bU} \cdot \bV \cdot \bV \cdot \bV = \frac{1}{2 \varepsilon^{3}}
    \bigg[ -\bF(\bU - 2 \varepsilon \bV) + 2\bF(\bU - \varepsilon \bV) - 2\bF(\bU + \varepsilon \bV)+ \bF(\bU + 2 \varepsilon \bV)
    \bigg] + \mathcal{O} ( \varepsilon^{2} ).
\end{equation}
We can further extend the procedure to compute the contraction
with three different vectors, \( \bV, \bW \), and \( \bX \),
\begin{equation}\label{eq:orig-sf-don-vec}
    \begin{split}
        \bF_{\bU \bU \bU} \cdot \bV \cdot \bW \cdot \bX = \frac{1}{6} \bigg[
            &\bF_{\bU \bU \bU} \cdot \left( \bV + \bW + \bX \right) \cdot \left( \bV + \bW + \bX \right) \cdot \left( \bV + \bW + \bX \right) \\
            & -\bF_{\bU \bU \bU} \cdot \left( \bV + \bW \right) \cdot \left( \bV + \bW \right) \cdot \left( \bV + \bW \right) \\
            & -\bF_{\bU \bU \bU} \cdot \left( \bV + \bX \right) \cdot \left( \bV + \bX \right) \cdot \left( \bV + \bX \right) \\
            & -\bF_{\bU \bU \bU} \cdot \left( \bW + \bX \right) \cdot \left( \bW + \bX \right) \cdot \left( \bW + \bX \right) \\
            & + \bF_{\bU \bU \bU} \cdot \bV \cdot \bV \cdot \bV +
                \bF_{\bU \bU \bU} \cdot \bW \cdot \bW \cdot \bW +
                \bF_{\bU \bU \bU} \cdot \bX \cdot \bX \cdot \bX
        \bigg],
    \end{split}
\end{equation}
and only to see that the number of terms to be computed 
rapidly increases in high-order tensor contraction terms.

As such, 
the original SF method in~\cref{eq:orig-sf-don-vec}
becomes less attractive for any PIF method higher than third-order accuracy, 
as it demands increasing
complexity in code implementation, which results in a significant loss
in the overall performance of the code.
For example, it requires 28 times flux function calls
for just getting a single tensor contraction,
\( \bF_{\bU \bU \bU} \cdot \bV \cdot \bW \cdot \bX \).
We should note that the major
bottleneck of the original system-free method stems from
\cref{eq:orig-sf-hes-vec,eq:orig-sf-don-vec}
that require to perform the \textit{Jacobian-like} approximations
multiple times.

To mitigation the computational bottleneck, in the following section, 
we introduce a newly improved version of
the system-free method, which does not require any further modifications,
even for the case of the tensor contractions of different vectors.

\section{A newly improved recursive system-free (SF) approach}\label{sec:new-sf}

In this chapter, we will improve the original system-free method~\cite{lee2021single}
to be applied in a recursive manner for higher-order derivatives of \( \bF \),
ensuing much simpler code complexity and faster performance.
Primarily, we define a functional \( \mathcal{D}_{u} \) that
represents the Jacobian-free method denoted in~\cref{eq:orig-sf-jac},
\begin{equation}\label{eq:jac-free}
    \bF_{\bU} \cdot \bV \approx \mathcal{D}_{u} (\bF \scolon \bV) \coloneqq
    \frac{1}{2\varepsilon_{v}} \bigg[
        \bF(\bU + \varepsilon_{v} \bV) - \bF(\bU - \varepsilon_{v} \bV)
    \bigg],
\end{equation}
where \( \varepsilon_{v} \) is the appropriately calculated \( \varepsilon \)
corresponding to the vector \( \bV \)
by following the original idea of the system-free method in~\cite{lee2021single},
\begin{equation}\label{eq:eps}
    \varepsilon_{v} = \min \left(\Delta t,\; \bar{\varepsilon}_{v} \right), \; \text{where} \;
    \bar{\varepsilon}_{v} = \frac{\sqrt{\varepsilon^{op}}}{\left\lVert \bV \right\rVert_{2}}.
\end{equation}
The study in this paper uses \( \varepsilon^{op} = \num{4.8062e-06} \) that is the optimal \( \varepsilon \)
value for the second-order Jacobian-free approximation in the 64-bit machine.
This choice is also justifiable for the recursive scheme considered below,
where 
the functional \( \mathcal{D}_{u} \) itself is defined as the Jacobian-free method fundamentally.
More detailed discussions about \( \varepsilon \) calculations
could be found in~\cite{lee2021single,an2011finite} and~\cite{knoll2004jacobian}.

We continue to apply \( \mathcal{D}_{u} \) in the following successive fashion to
calculate the tensor contractions between higher-order derivatives for
the flux function \( \bF \) and arbitrary vectors.
Thus, the Hessian approximation is,
\begin{equation}\label{eq:hes-free}
    \begin{split}
        \bF_{\bU \bU} \cdot \bV \cdot \bW &\approx
        \mathcal{D}_{u} \Big( \mathcal{D}_{u} (\bF \scolon \bV) \scolon \bW \Big) \\
        &
        \begin{split}
            =\frac{1}{4 \varepsilon_{v} \varepsilon_{w}}
                \bigg[
                     &\bF(\bU + \varepsilon_{v} \bV + \varepsilon_{w} \bW)
                    -\bF(\bU - \varepsilon_{v} \bV + \varepsilon_{w} \bW)\\
                    -&\bF(\bU + \varepsilon_{v} \bV - \varepsilon_{w} \bW)
                    +\bF(\bU - \varepsilon_{v} \bV - \varepsilon_{w} \bW)
                \bigg].
        \end{split}
    \end{split}
\end{equation}
Again, following~\cref{eq:eps},
\( \varepsilon_{v} \) and \( \varepsilon_{w} \) are the optimal \( \varepsilon \) values
normalized by its corresponding vectors \( \bV \) and \( \bW \), respectively.

Note that the improved version of the Hessian-free method in \cref{eq:hes-free}
is applicable regardless the tensor contraction is
with two identical vectors (e.g., \( \bF_{\bU \bU} \cdot \bV \cdot \bV \))
or with two distinct vectors (e.g., \( \bF_{\bU \bU} \cdot \bV \cdot \bW \)),
hence it does not require separate formulations as in
\cref{eq:orig-sf-hes-vec,eq:orig-sf-don-vec}.

The simplicity gain from the improved version of the system-free method
is further rewarded when we apply it to
the higher-order derivatives of \( \bF \).
Following the equivalent strategy, the tensor contraction of
the third-order derivative of the flux function, \( \bF_{\bU \bU \bU} \)
with three distinct vectors, \( \bV, \bW \), and \( \bX \) is written compactly as,
\begin{equation}\label{eq:don-free}
    \begin{split}
        \bF_{\bU \bU \bU} \cdot \bV \cdot \bW \cdot \bX &\approx
        \mathcal{D}_{u} \Bigg( \mathcal{D}_{u} \Big( \mathcal{D}_{u} (\bF \scolon \bV) \scolon \bW \Big) \scolon \bX \Bigg) \\
        &
        \begin{split}
            =\frac{1}{8 \varepsilon_{v} \varepsilon_{w} \varepsilon_{x}}
            \bigg[
                &\bF(\bU + \varepsilon_{v} \bV + \varepsilon_{w} \bW + \varepsilon_{x} \bX)
                -\bF(\bU - \varepsilon_{v} \bV + \varepsilon_{w} \bW + \varepsilon_{x} \bX)\\
                -&\bF(\bU + \varepsilon_{v} \bV - \varepsilon_{w} \bW + \varepsilon_{x} \bX)
                +\bF(\bU - \varepsilon_{v} \bV - \varepsilon_{w} \bW + \varepsilon_{x} \bX)\\
                -&\bF(\bU + \varepsilon_{v} \bV + \varepsilon_{w} \bW - \varepsilon_{x} \bX)
                +\bF(\bU - \varepsilon_{v} \bV + \varepsilon_{w} \bW - \varepsilon_{x} \bX)\\
                +&\bF(\bU + \varepsilon_{v} \bV - \varepsilon_{w} \bW - \varepsilon_{x} \bX)
                -\bF(\bU - \varepsilon_{v} \bV - \varepsilon_{w} \bW - \varepsilon_{x} \bX)
            \bigg].
        \end{split}
    \end{split}
\end{equation}

Let us consider the number of flux function calls of SF approximations
to compare the performance gain from the recursive SF method.
Comparing~\cref{eq:orig-sf-don-vec} and~\eqref{eq:don-free}, for example,
the numerical flux function needed to be called 28 times
in the original SF method. However, the
recursive SF method only requires eight evaluations.
This is a huge improvement in both performance and compactness. 

With the modified SF method in~\crefrange{eq:jac-free}{eq:don-free},
we can approximate all the tensor contractions
needed for calculating temporal flux derivatives
in~\crefrange{eq:Ft}{eq:divtx} without analytical
derivations for \textit{Jacobian-like} terms,
giving the system independence of the high-order scheme.
We should note that the recursive modifications of the SF method
presented in this section do not affect the
solution's accuracy and stability compared to the original
SF method.

\section{Numerical Test Results}\label{sec:results}

This section presents numerical results of
well-known test problems for benchmarking the modified SF-PIF
methods' capabilities. We mainly compare the results
from the third-order and the fourth-order temporal methods of SSP-RK and
SF-PIF. We used the well-known three-stages, third-order
SSP-RK method~\cite{gottlieb1998total} and
the five-stages, fourth-order SSP-RK method~\cite{spiteri2002new}
(RK3 and RK4 hereafter)
for the comparisons of the third-order and fourth-order SF-PIF methods
(SF-PIF3 and SF-PIF4 hereafter),
respectively.

The main purpose of this section is to demonstrate that the proposed
recursive SF-PIF methods produce the same quality of solutions,
with less CPU time in the light of the single-stage time update,
as compared to
RK3 and RK4 methods applied to FDM discretizations.
We use the conventional fifth-order WENO-JS~\cite{jiang1996efficient}
spatial reconstruction method for all simulations; thus, we expect
fifth-order spatial accuracy \( \mathcal{O}(\dx^{5}) \) combined with
third \( \mathcal{O}(\dt^{3}) \) or fourth-order \( \mathcal{O}(\dt^{4}) \) temporal accuracy
for all numerical results.
The fixed Courant numbers, \( C_{\text{cfl}} = 0.4 \) and \( C_{\text{cfl}} = 0.3 \),
are used for all 2D and 3D simulations, respectively.

\subsection{2D Euler equations}\label{sec:results-2d}

\subsubsection{The Sod's shock tube test (rotated \(\ang{45}\))}\label{sec:sod45}

\begin{figure}
    \centering
    \begin{subfigure}{80mm}
        \centering
        \includegraphics[width=0.95\textwidth]{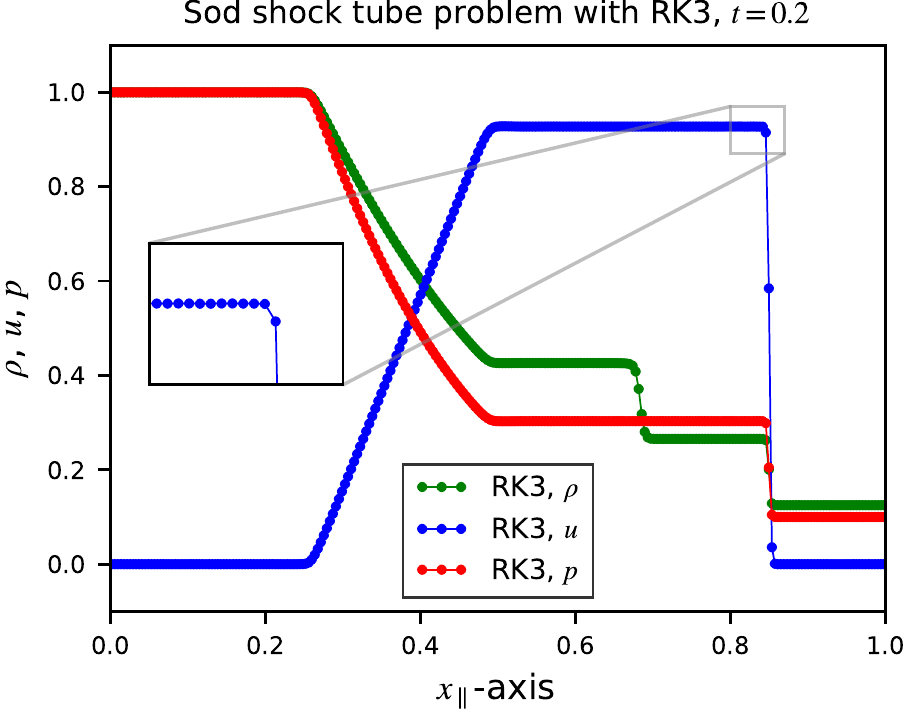}
    \end{subfigure}
    \begin{subfigure}{80mm}
        \centering
        \includegraphics[width=0.95\textwidth]{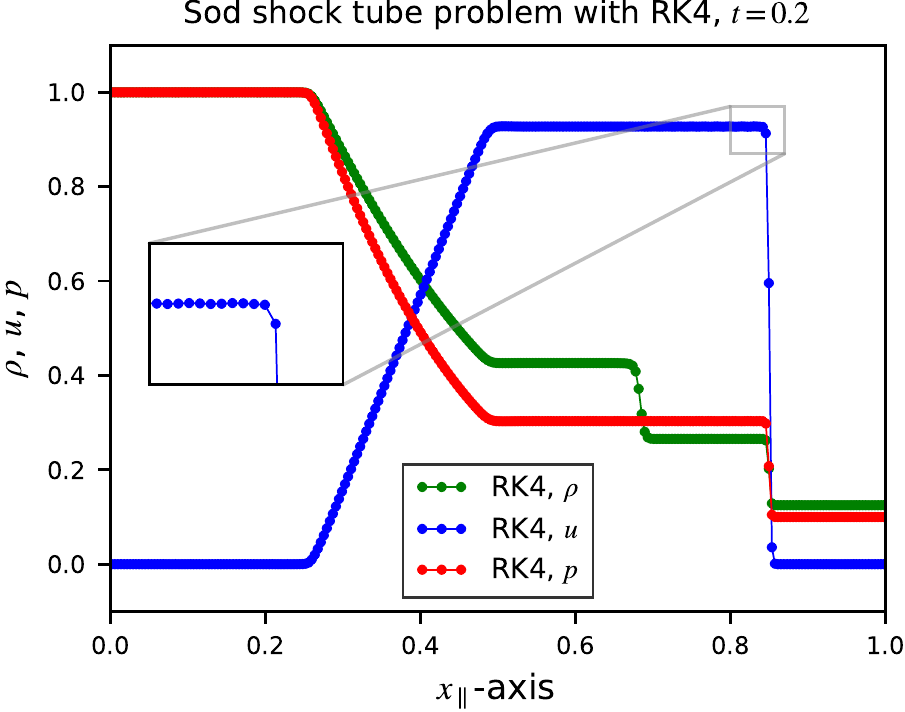}
    \end{subfigure} 
    \begin{subfigure}{80mm}
        \centering
        \includegraphics[width=0.95\textwidth]{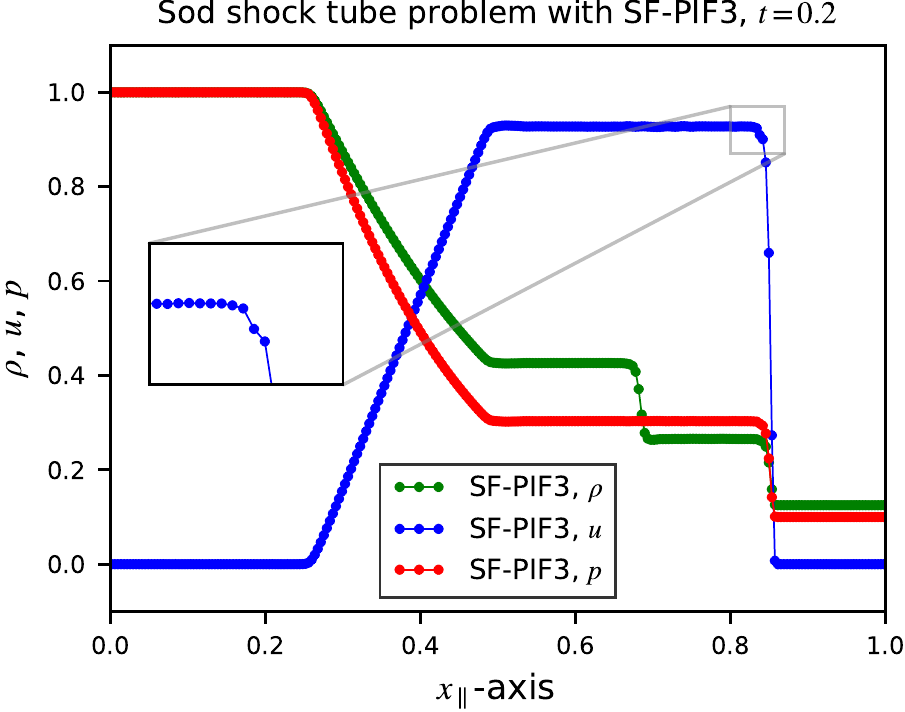}
    \end{subfigure}
    \begin{subfigure}{80mm}
        \centering
        \includegraphics[width=0.95\textwidth]{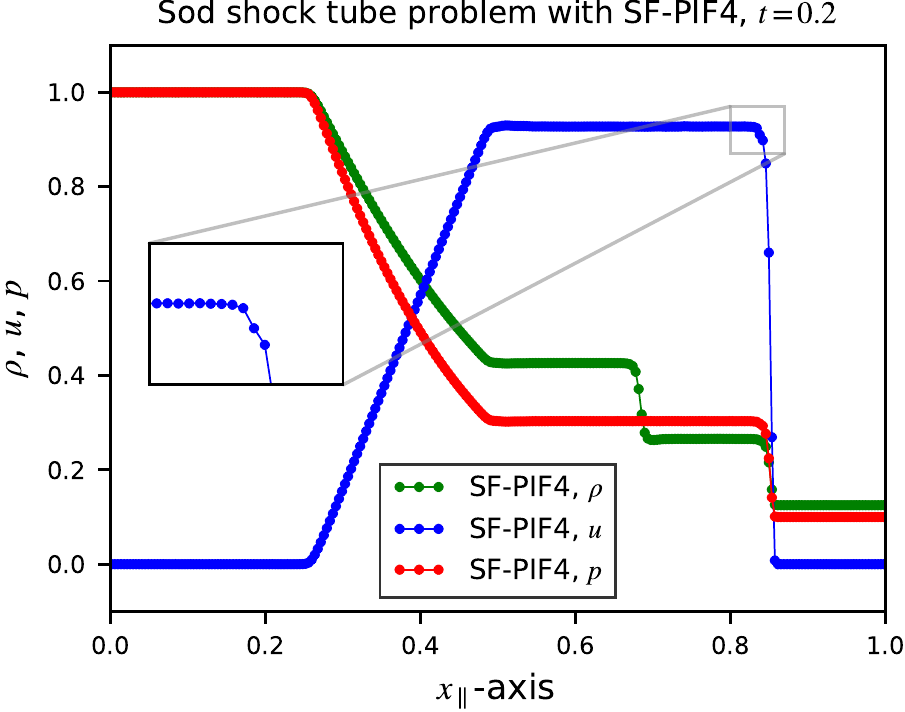}
    \end{subfigure}
    \caption{One-dimensional profiles
        of the inclined Sod's shock tube problem along the \( x_{\parallel} \) direction,
        displaying the fluid density $\rho$, the $x$-velocity $u$, and the thermal pressure $p$,
        illustrated at \( t=0.2 \).
        All simulations performed in a 2D simulation box of
        \( [1024 \times 1024] \) grid resolution,
        but the number of data points on each panel is \( N_{\parallel} \) = 256,
        as we only take the first quarter of the diagonal axis for the shown profiles.
        \textbf{Left column}: The profiles solved by RK3 (top) and SF-PIF3 (bottom).
        \textbf{Right column}: The profiles solved by RK4 (top) and SF-PIF4 (bottom).
    }\label{fig:sod45}
\end{figure}

We start with the classical shock tube problem of Sod \cite{sod1978survey}, testing a numerical scheme's ability to capture the three wave families of the 1D Riemann problem: a shock, contact discontinuity, and rarefaction fan.
Although Sod's problem is a 1D shock tube problem originally,
we implement the test in a 2D domain by tilting the shock wave
direction by the angle of \( \theta = \ang{45} \). We adopt the idea of
Kawai~\cite{kawai2013divergence}, where the initial conditions are
repeated multiple times along the direction of the wave propagation 
so that the problem may be executed with periodic boundary conditions.
Explicitly, the initial condition is given as,
\begin{equation}\label{eq:sod45}
    \left( \rho, u, v, p \right) = \begin{cases}
        \left( 1, 0, 0, 1 \right) & \text{for } x_{\parallel} \le 0.5, \quad 1.5 < x_{\parallel} \le 2.5, \quad 3.5 < x_{\parallel} \le 4, \\
        \left( 0.125, 0, 0, 0.1 \right) & \text{for } 0.5 < x_{\parallel} \le 1.5, \quad 2.5 < x_{\parallel} \le 3.5,
    \end{cases}
\end{equation}
where \( x_{\parallel} = x \cos{\theta} + y \sin{\theta} \) is the direction
parallel to the wave propagation. The simulation domain is a periodic box of
\( [0, 2/\cos{\theta}] \times [0, 2/\sin{\theta}] \).
For the final result profiles, we take only the
bottom-left quarter of the diagonal axis, \( 0 \le x_{\parallel} \le 1 \);
thus, the number of data points of the result profiles
would be a quarter of the grid resolution in \( x \) and \( y \),
\( N_{\parallel} = N_{x}/4 = N_{y}/4 \).

The \( \ang{45} \)-angled Sod's shock tube test results at \( t = 0.2 \)
are depicted in~\cref{fig:sod45}. We used the grid resolution of \( N_{x} = N_{y} = 1024 \);
thus, the figure shows an \( N_{\parallel} = 256 \) number of data points
along with the diagonal axis, \( x_{\parallel} \). As shown in~\cref{fig:sod45},
SF-PIF methods produce fairly comparable results to those of RK3 and RK4,
except for the small oscillation in the shock front of \( x \)-velocity.
This issue was already discussed in our previous study~\cite{lee2021single},
where the oscillations are originated from the central differencing
formulae that we used for getting spatial derivatives. We can minimize
the oscillation by applying WENO-\textit{like} differencing
introduced in the appendix of~\cite{lee2021single}. However,
we choose to use the conventional central differencing formulae
for this study as we did not observe any unphysical oscillations
for the rest of the test problems.

\subsubsection{The Shu-Osher problem (rotated \(\ang{45}\))}\label{sec:shu45}

Our next choice of test problem is the Shu-Osher problem~\cite{shu1989efficient}
that describes the interactions between a Mach 3 shock and a smooth density profile.
Initially, a Mach 3 shock wave travels to the right through a sinusoidally perturbed
density profile. As the shock wave propagates along the perturbed region,
the profile gets compressed, resulting in a frequency-doubled region behind the shock. 
As the shock wave moves further to the right,
the doubled-frequency region returns to its original frequency, at which point
it becomes a sequence of sharp profile instead of smooth sine wave
due to the shock-steepening.

We performed the Shu-Osher problem in 2D by inclining the shock wave direction
by an angle of \(\theta = \ang{45} \), similar to what we did in the~\cref{sec:sod45}.
We use periodic boundary conditions in both directions
of the simulation domain \( [0, 20/\cos{\theta}] \times [0, 20/\sin{\theta}] \).
The grid resolution is \( N_{x} = N_{y} = 1024 \); thus, the effective
resolution \( N_{\parallel} \) is 256 as we take
only the bottom-left quarter of \( x_{\parallel} \) to report the 1D profiles.

\begin{figure}
    \centering
    \begin{subfigure}{80mm}
        \centering
        \includegraphics[width=0.95\textwidth]{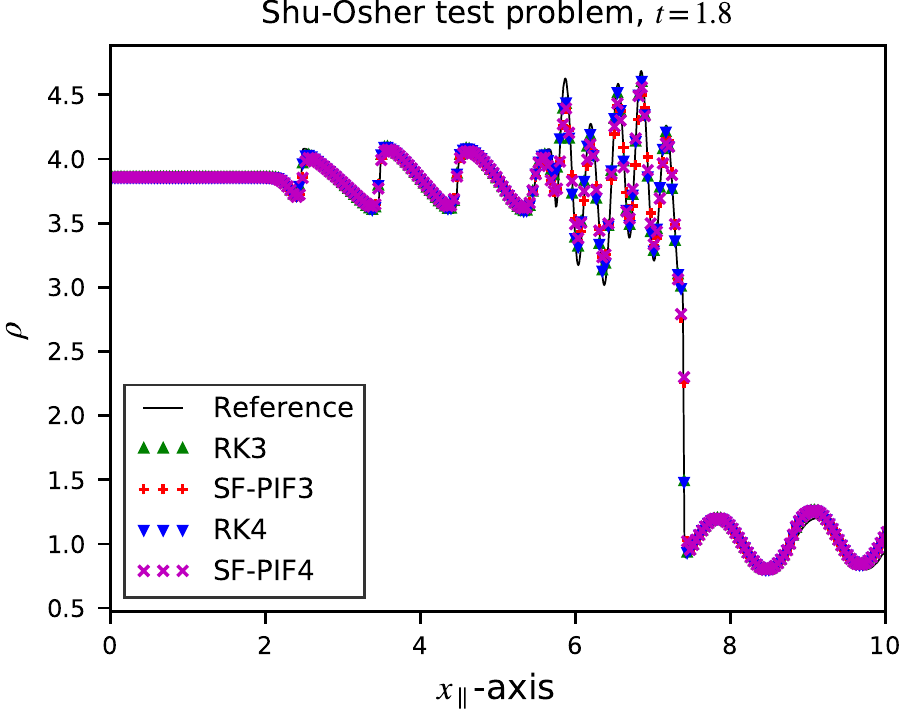}
    \end{subfigure}
    \begin{subfigure}{80mm}
        \centering
        \includegraphics[width=0.95\textwidth]{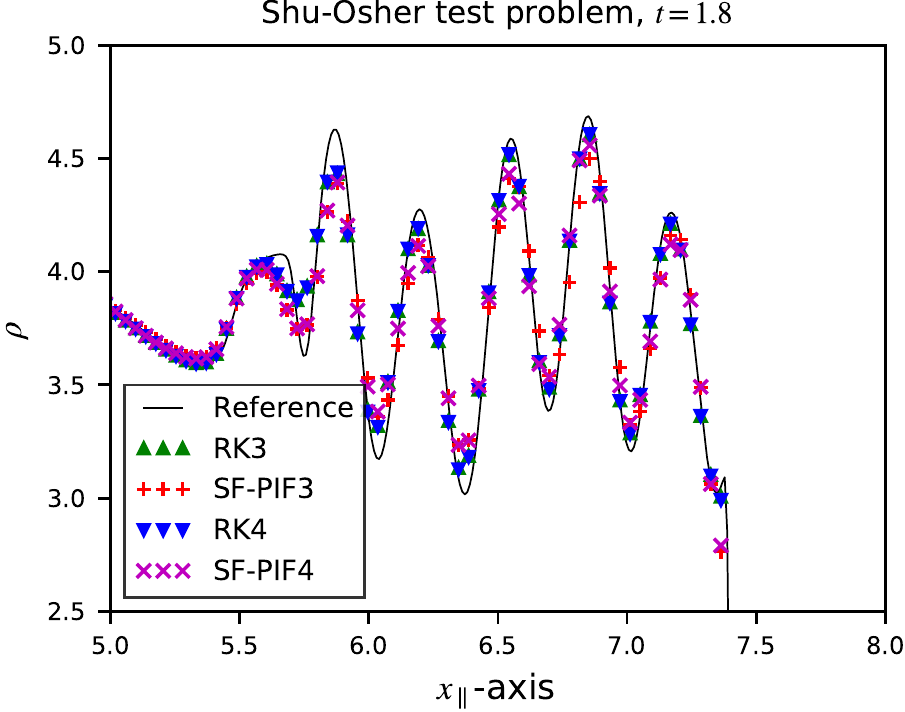}
    \end{subfigure}
    \caption{One dimensional density profiles along the \( x_{\parallel} \) direction 
        of the inclined Shu-Osher problem at \( t = 1.8 \).
        The solid line represents the reference solution,
        solved by RK4 with 1024 data points
        in the diagonal axis, i.e., \( N_{\parallel} = 1024 \).
        All other solutions, represented by the symbols, are resolved on an
        \( N_{\parallel} = 256 \) grid resolution in the diagonal axis.
        The detailed view of the high-frequency region is shown
        on the right panel.
    }\label{fig:shu45}
\end{figure}

The density profiles along the \( x_{\parallel} \) direction
are given in~\cref{fig:shu45}.
The four different temporal method choices (RK3, RK4, SF-PIF3, and SF-PIF4)
produce reasonably acceptable solution profiles capturing the high-frequency 
 amplitudes fairly well in the frequency-doubled region. We see that the results of RK3 and RK4 are nearly
identical, while SF-PIF3 and SF-PIF4 show slight differences near the highest amplitudes.
Generally, RK methods capture the amplitudes marginally better,
but in the left-most part of the double-frequency region, \( x \approx 5.8 \),
SF-PIF methods capture the highest peak of the amplitude better than
the RK methods near the transition between the frequency 
doubling and the shock steepened perturbations.

\subsubsection{Isentropic vortex advection}\label{sec:vortex}

The isentropic vortex advection problem~\cite{shu1998essentially}
is one of the most popular test choices to measure the simulation code's
accuracy and performance. Although the problem is fully nonlinear, 
the exact solution is always existent in the form of its initial condition, from which
an isentropic vortex is advected through periodic boundaries. 
We can evaluate the accuracy of a method on nonlinear problems by comparing
the final result with its initial condition.
Here, we adopt the idea from Spiegel~\cite{spiegel2015survey}, where the
simulation domain size is doubled up as \( [0, 20] \times [0, 20] \)
to prevent vortex-vortex couplings near the periodic boundaries.

\begin{figure}
    \centering
    \begin{subfigure}{80mm}
        \centering
        \includegraphics[width=0.95\textwidth]{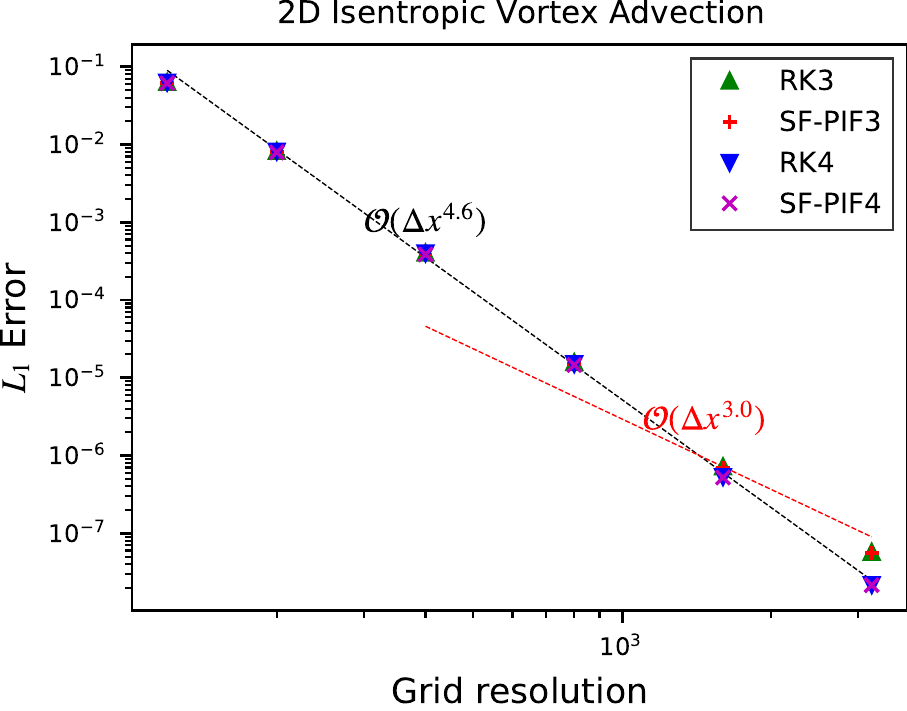}
    \end{subfigure}
    \begin{subfigure}{80mm}
        \centering
        \includegraphics[width=0.95\textwidth]{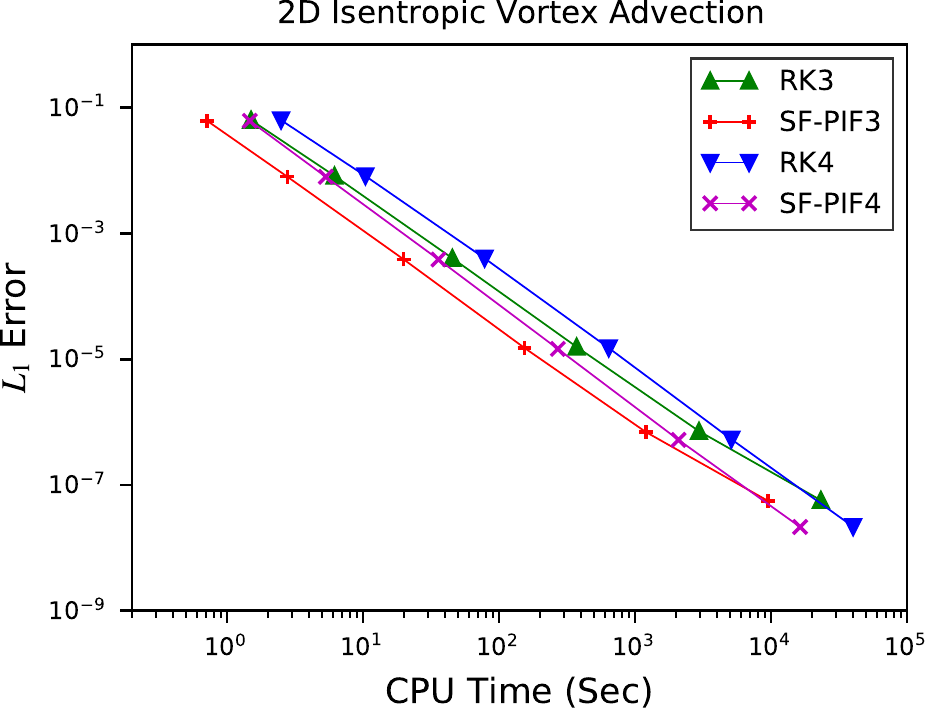}
    \end{subfigure}
    \caption{The \( L_{1} \) errors of the isentropic vortex advection test problem
        with respect to the grid resolutions (\textbf{left});
        with respect to the computation time (\textbf{right}).
    }\label{fig:vortex}
\end{figure}

The results of \( L_{1} \) errors on six different grid resolutions,
\( N_{x} = N_{y} = 120, 200, 400, 800, 1600, \) and \( 3200 \), are plotted
in the left panel of~\cref{fig:vortex}. As expected, all temporal methods
follow the convergence line of order \( \sim \mathcal{O}(\dx^{4.6}) \),
which is nearly the same as WENO's fifth-order spatial accuracy.
However, at the critical grid resolution, \( N_{x} = N_{y} = 1600 \),
the third-order temporal methods of RK3 and SF-PIF3 
start to degrade the whole solution accuracy.
This behavior can be explained that the errors from the fifth-order spatial WENO
solver are dominant on the grid resolutions 
up to  \( N_{x} = N_{y} = 1600 \), after which
the truncation errors associated with the third-order temporal integrators 
become dominant over
%
%
%
the error of the fifth-order spatial solver.
This result tells us the significance of high-order temporal updates 
in fine grid resolutions: 
a high-order spatial method does
require a \textit{comparably} high-order temporal method
to retain the overall quality of the solutions,
particularly when we are motivated to add more grid 
resolutions to resolve finer scales more accurately.
Otherwise, the temporal solver's lower order accuracy 
can potentially degrade the solution accuracy, 
contradicting the intended motivation.
%
A direct consequence of this observation applies to
CFD simulations on an adaptive mesh refinement (AMR) configuration.
A mediocre low-order temporal method is integrated with a high-order spatial solver.
This combination of such two solvers
creates a computational dilemma of not making any further enhancement
in solution accuracy 
as the computational grids are progressively refined to 
improve the quality of the AMR solutions. Instead, the solution 
accuracy is to be bounded by the lower temporal accuracy.

%
%

To compare the computational performance over the four different
temporal integrators, we present \( L_{1} \) errors as a function of CPU time data
on the right panel of~\cref{fig:vortex}. 
The drop of the convergence rates is again observed 
in the two third-order solutions of RK3 and SF-PIF3 
at the critical resolution \( N_{x} = N_{y} = 1600 \).
On the contrary, the two fourth-order solutions of RK4 
and SF-PIF4 continue at fourth-order without any change.
Up to this resolution,
SF-PIF3 is the fastest in reaching any given target \( L_{1} \) error threshold
in CPU time.
However, 
on any grid resolutions finer than the critical resolution,
SF-PIF3's \( L_{1} \) error drops to the third-order convergence, which ultimately
crosses the SF-PIF4's fourth-order straight convergence line at
\( N_{x} = N_{y} = 10^4 \), beyond which
its error will remain larger than the ones from the fourth-order
temporal schemes as long as the convergence rate follows the
pattern at the high-resolution tail.
\begin{table}[hb!]
   \footnotesize
    \centering
    \caption{The \( L_{1} \) errors, the rates of convergence,
        and the computation times for the vortex advection test
        solved using RK3 and SF-PIF3 methods (\textbf{top});
        using RK4 and SF-PIF4 methods (\textbf{bottom}).
        All simulation runs are performed on the four 20-cores
        Cascade Lake Intel Xeon processors, utilized 64 parallel threads.
        CPU times are measured in seconds, averaged over 10 individual runs.
    }\label{table:vortex3-4}
    \begin{tabular}{@{}ccccclcccc@{}}
        \toprule
        \multirow{2}{*}{\( N_{x} = N_{y} \)} & \multicolumn{4}{c}{RK3} &  & \multicolumn{4}{c}{SF-PIF3} \\
        \cmidrule(lr){2-5} \cmidrule(l){7-10}
        & \(L_{1}\) error & \(L_{1}\) order & CPU Time & Speedup &  &
        \(L_{1}\) error & \(L_{1}\) order & CPU Time & Speedup \\ \midrule
        120  & \num{6.31E-2} & \--- & \SI{1.50}{\second}      & 1.0 &  & \num{6.16E-2} & \--- & \SI{0.71}{\second}    & 0.48 \\
        200  & \num{8.20E-3} & 4.00 & \SI{6.17}{\second}      & 1.0 &  & \num{7.96E-3} & 4.00 & \SI{2.77}{\second}    & 0.45 \\
        400  & \num{4.02E-4} & 4.35 & \SI{45.44}{\second}     & 1.0 &  & \num{3.86E-4} & 4.37 & \SI{19.89}{\second}   & 0.44 \\
        800  & \num{1.57E-5} & 4.68 & \SI{372.47}{\second}    & 1.0 &  & \num{1.51E-5} & 4.68 & \SI{153.92}{\second}  & 0.41 \\
        1600 & \num{7.18E-7} & 4.45 & \SI{2957.26}{\second}   & 1.0 &  & \num{6.95E-7} & 4.44 & \SI{1203.10}{\second} & 0.41 \\
        3200 & \num{5.72E-8} & 3.65 & \SI{23274.37}{\second}  & 1.0 &  & \num{5.60E-8} & 3.63 & \SI{9494.65}{\second} & 0.41 \\
    \end{tabular}
    \begin{tabular}{@{}ccccclcccc@{}}
        \toprule
        \multirow{2}{*}{\( N_{x} = N_{y} \)} & \multicolumn{4}{c}{RK4} &  & \multicolumn{4}{c}{SF-PIF4} \\
        \cmidrule(lr){2-5} \cmidrule(l){7-10}
        & \(L_{1}\) error & \(L_{1}\) order & CPU Time & Speedup &  &
        \(L_{1}\) error & \(L_{1}\) order & CPU Time & Speedup \\ \midrule
        120  & \num{6.30E-2} & \--- & \SI{2.50}{\second}       & 1.0 &  & \num{6.14E-2} & \--- & \SI{1.47}{\second}     & 0.59 \\
        200  & \num{8.15E-3} & 4.00 & \SI{10.42}{\second}      & 1.0 &  & \num{7.91E-3} & 4.01 & \SI{5.33}{\second}     & 0.51 \\
        400  & \num{4.01E-4} & 4.35 & \SI{78.47}{\second}      & 1.0 &  & \num{3.85E-4} & 4.36 & \SI{35.89}{\second}    & 0.46 \\
        800  & \num{1.51E-5} & 4.73 & \SI{641.50}{\second}     & 1.0 &  & \num{1.46E-5} & 4.72 & \SI{270.94}{\second}   & 0.42 \\
        1600 & \num{5.33E-7} & 4.82 & \SI{5115.47}{\second}    & 1.0 &  & \num{5.21E-7} & 4.81 & \SI{2091.20}{\second}  & 0.41 \\
        3200 & \num{2.17E-8} & 4.62 & \SI{40195.034}{\second}  & 1.0 &  & \num{2.15E-8} & 4.60 & \SI{16377.73}{\second} & 0.41 \\
    \end{tabular}
\end{table}

On the other hand,
it is distinctively superior to see that
SF-PIF4's solution reaches any fixed target error
in a \textit{faster} CPU time than the
\textit{third-order} RK3's solution while keeping the
numerical errors as low as RK4 results at each grid resolution. The detailed
simulation data from the vortex advection tests are
presented in ~\cref{table:vortex3-4}.
All performance results are measured in second,
averaged over 10 simulation runs conducted on
two nodes of the UC Santa Cruz's high-performance computer, \textit{lux}.
Each node has two 20-core Cascade Lake Intel Xeon processors,
and we utilized 64 parallel threads for each simulation run.

\subsubsection{2D Riemann problem: Configuration 3}\label{sec:2drp}
To test the methods' ability to capture the complex
fluid structures in 2D, we performed a two-dimensional Riemann problem
called Configuration 3.
This specific setup and other types of configurations have been
extensively studied in~\cite{zhang1990conjecture, schulz1993classification, schulz1993numerical}
and have been adopted as popular benchmark test problems.
To set up Configuration 3, we follow the initial condition
from~\cite{don2016hybrid, lee2017piecewise,lee2021single}. 
We conducted numerical experiments on
a \( 1600 \times 1600 \) grid resolution, which is generally considered as 
a very high resolution in 2D. This grid resolution choice is made based on
our observation in~\cref{sec:vortex} 
to make sure that the temporal errors are comparable to
or dominant over the spatial errors; 
thereby, we can anticipate temporal error dominant results
that allow us to focus on the effects of the four different temporal solvers.

\begin{figure}[ht!]
    \centering
    \begin{subfigure}{160mm}
        \centering
        \includegraphics[width=0.9\textwidth]{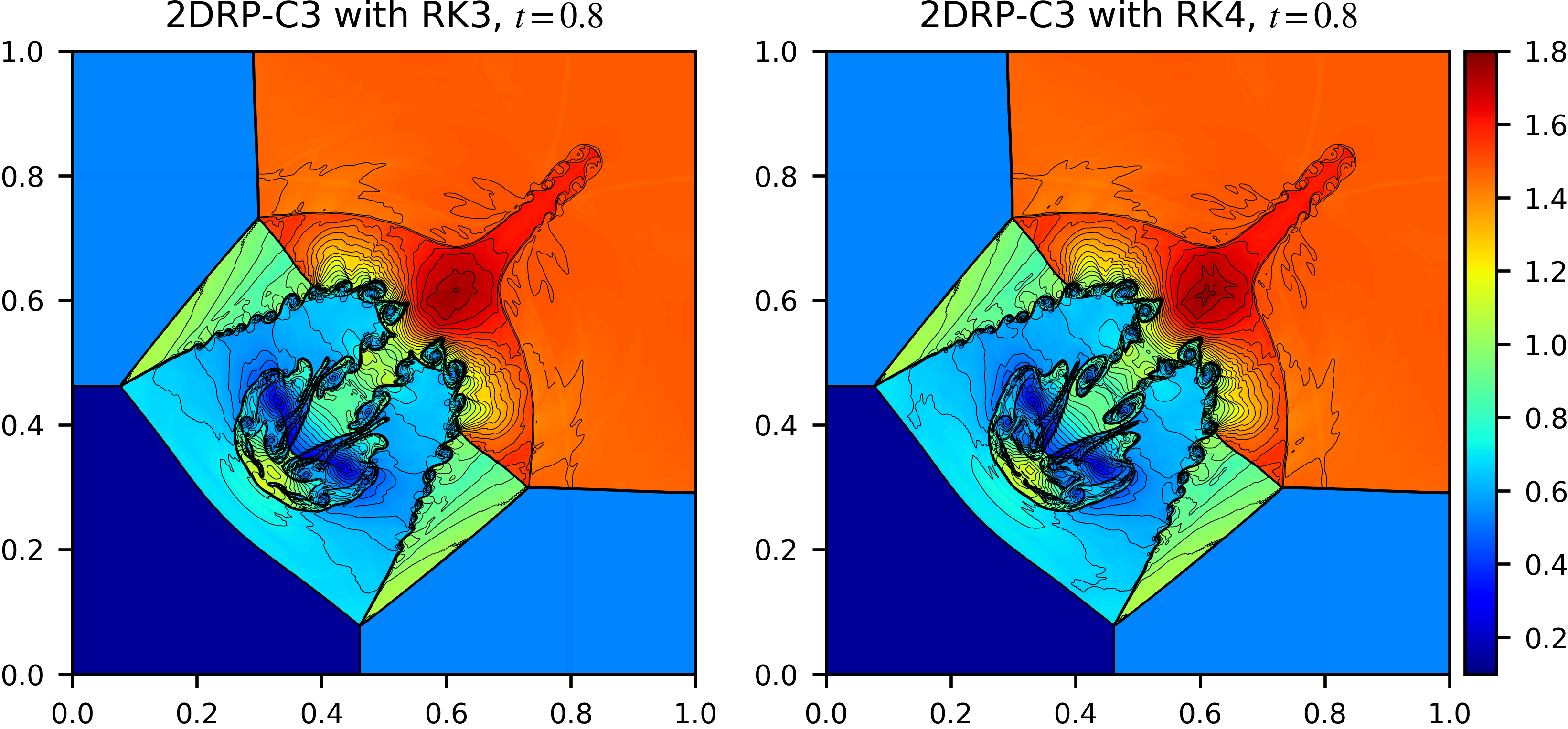}
    \end{subfigure}
    \begin{subfigure}{160mm}
        \centering
        \includegraphics[width=0.9\textwidth]{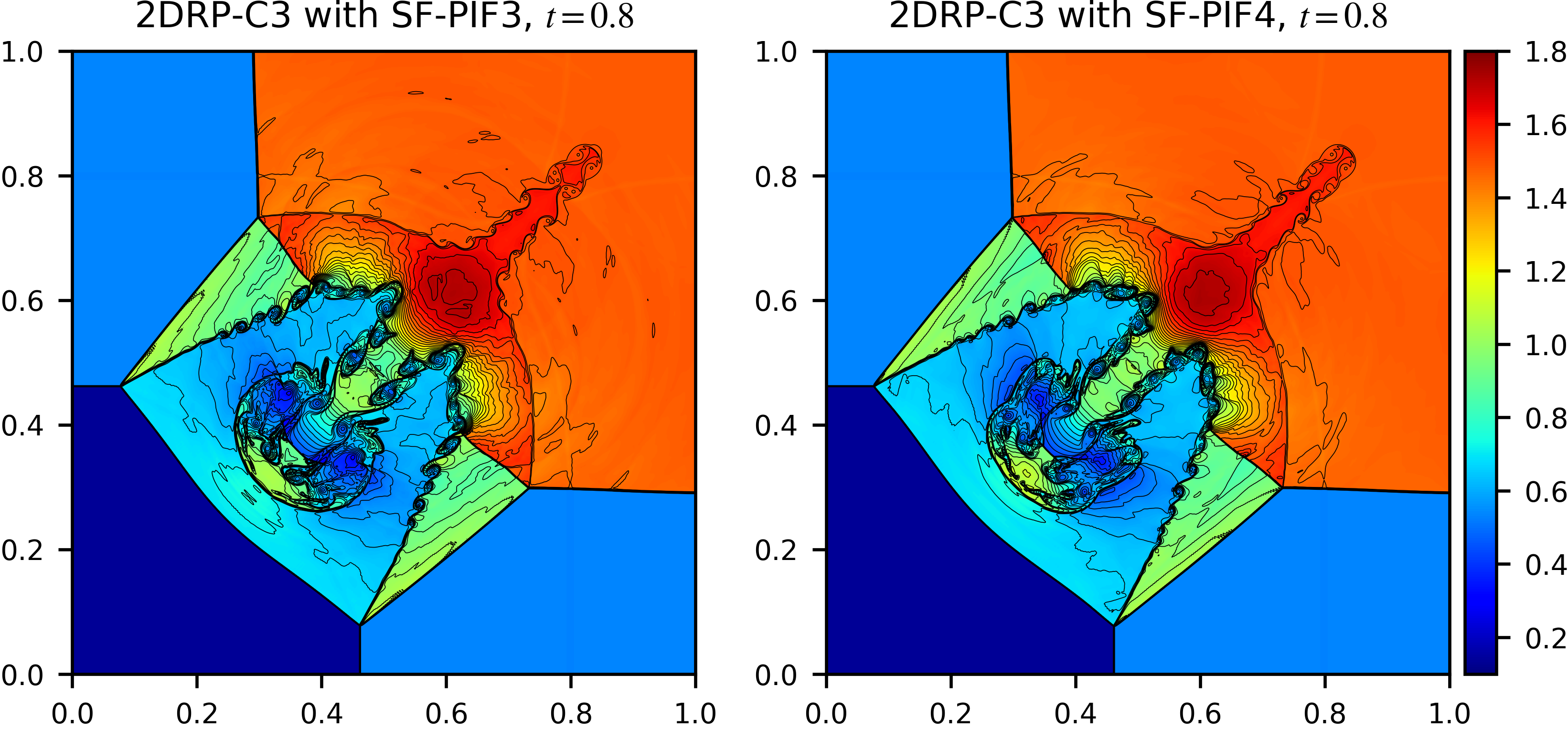}
    \end{subfigure}
    \caption{The density maps of Configuration 3 at \( t = 0.8 \).
        \textbf{Left column:} The solutions using RK3 (top) and SF-PIF3 (bottom).
        \textbf{Right column:} The solutions using RK4 (top) and SF-PIF4 (bottom).
        Forty levels of black contour lines are over-plotted in each figure
        with the same range of the color map.
        All simulations are performed on a \( 1600 \times 1600 \) grid resolution.
        }\label{fig:2drp_c3}
\end{figure}

The results at \( t = 0.8 \) are shown in~\cref{fig:2drp_c3}.
The pseudo-colors represent the density map ranging between \( [0.1, 1.8] \), and
40 contour lines within the same range are over-plotted as solid black lines.
We see that all four different temporal schemes produce well-known, acceptable results,
keeping the assumed diagonal symmetry exceptionally well on this high resolution.
This problem is highly nonlinear, involving formations of the upward-moving jet,
the downward-moving mushroom-jet,
secondary Kelvin-Helmholtz instabilities exhibited
as the small-scale vortical rollups along the slip lines and along the stems of the two jets.
As such, it is a non-trivial task to address if
a method under consideration is \textit{better} or \textit{worse} based on
the number of such rollups in the simulations (see our discussion 
in~\cite{lee2021single}). At best, such quantification can only provide
proof of intrinsic information about the amount of 
numerical dissipation of each method.
From this perspective, we conclude that the two SF-PIF solutions
produce the equivalent amount of vortical rollups compared with the corresponding
RK solutions, confirming the SF-PIF methods' validity.

\subsubsection{Double Mach reflection}\label{sec:dmr}
Our final 2D test case is the double Mach reflection test
that launches a strong Mach 10 shock with \(\ang{60}\) of an incident angle
between the shock plane and the bottom reflecting wedge.
The initial condition is the same as the original setup introduced by Woodward and Colella~\cite{woodward1984numerical},
except that we doubled the \( y \)-domain size following~\cite{kemm2016proper}
to prevent numerical artifacts from the top boundary interaction with the
secondary shock wave and the slip line.

\begin{figure}
    \centering
    \begin{subfigure}{160mm}
        \centering
        \includegraphics[width=0.95\textwidth]{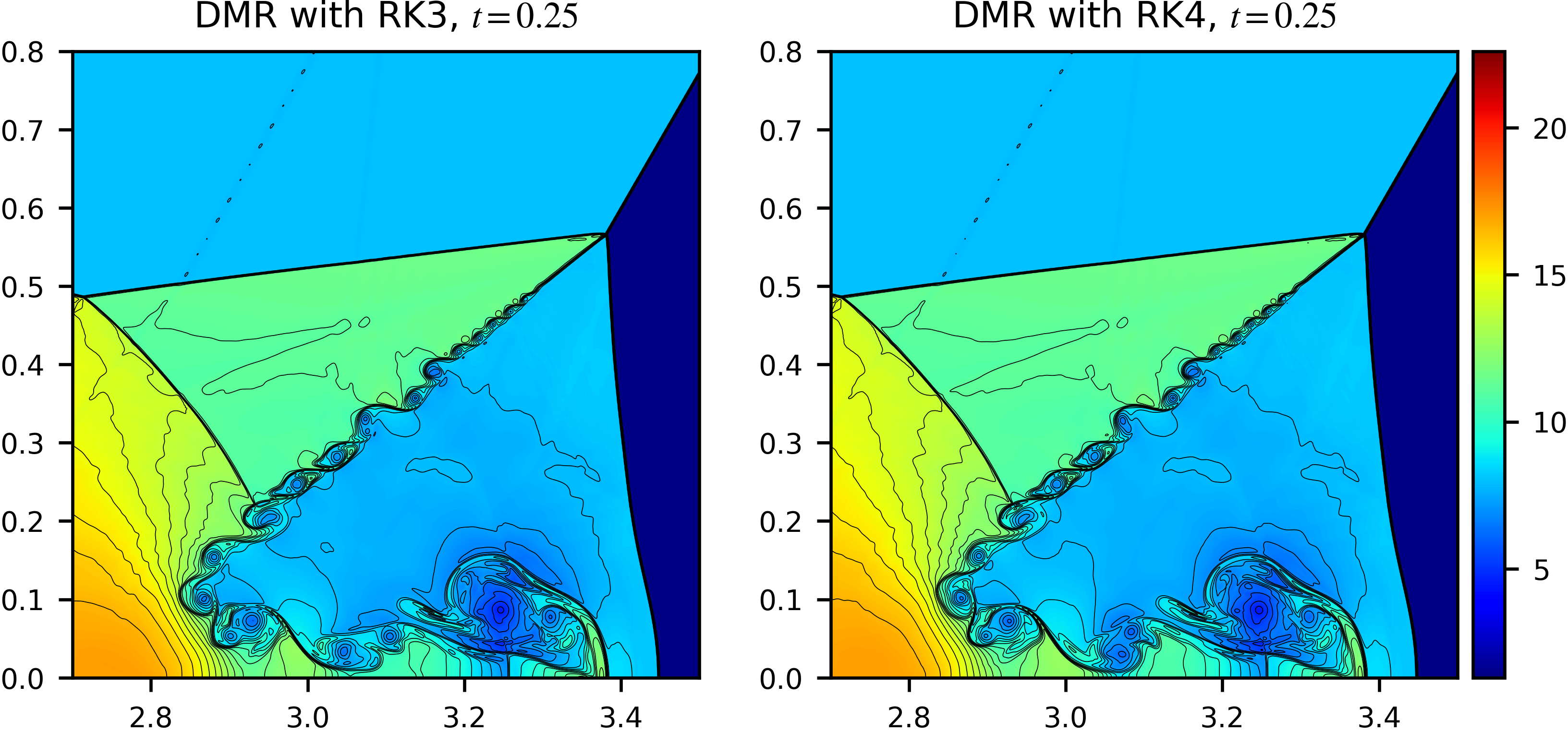}
    \end{subfigure}
    \begin{subfigure}{160mm}
        \centering
        \includegraphics[width=0.95\textwidth]{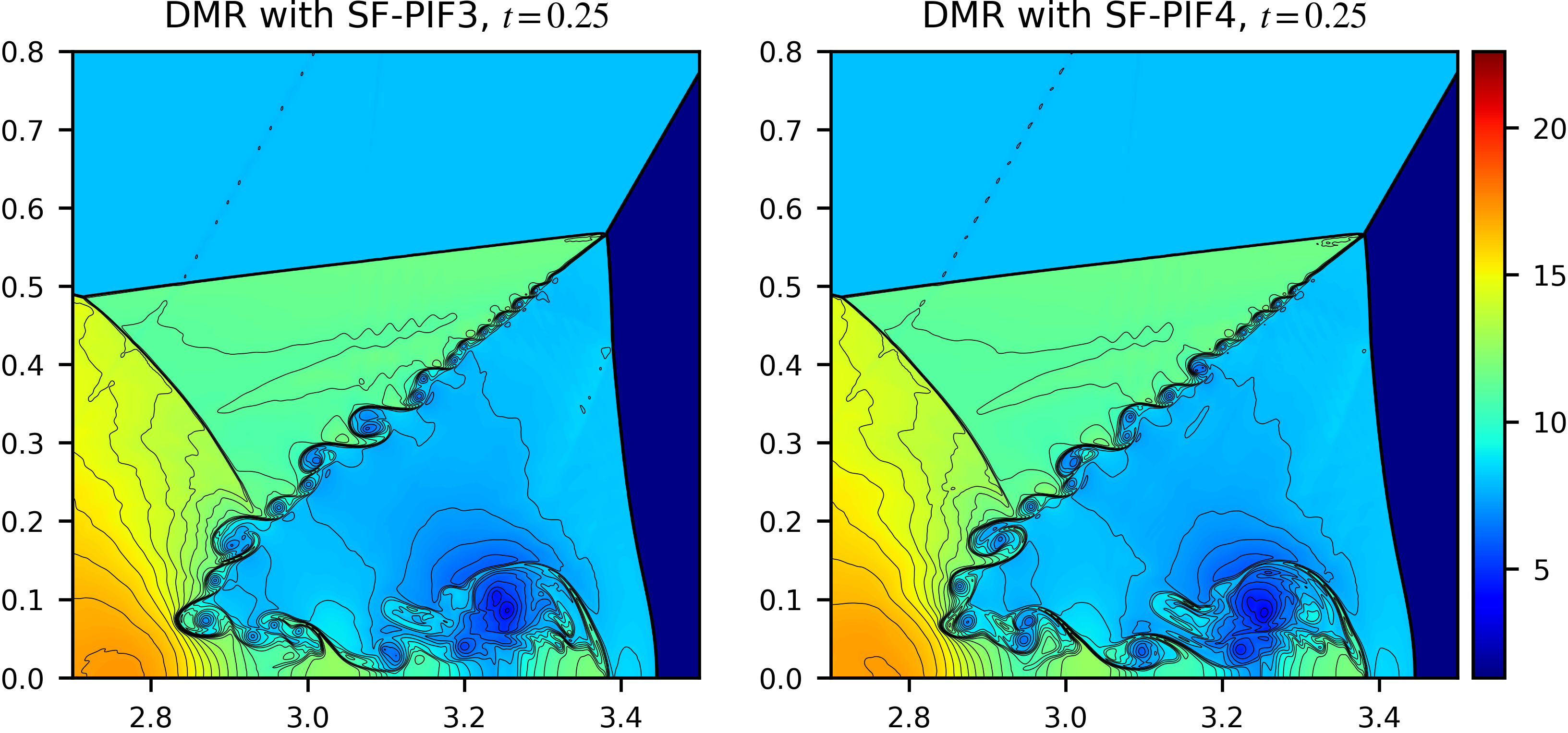}
    \end{subfigure}
    \caption{The density map of the double Mach reflection test at \( t = 0.25 \)
        zoomed-in near the jet. Forty levels of contour lines are over-plotted
        in solid black curves with the same range of the color map.
        All simulation results are performed on a
        \( 4096 \times 2048 \) grid resolution.
        \textbf{Left column:} The solutions using RK3 (top) and SF-PIF3 (bottom).
        \textbf{Right column:} The solutions using  RK4 (top) and SF-PIF4 (bottom).
    }\label{fig:dmr}
\end{figure}

The density results are presented in~\cref{fig:dmr}. The density color map
ranges between \( [1.3, 22.6] \), and the solid black lines represent
40 levels of the density contour lines with the same range.
The figures are zoomed-in near the vicinity of the jet and the primary
triple point, which is widely considered the main area of interest in
the Double Mach Reflection test. All simulation runs are performed on a
\( 4096 \times 2048 \) grid resolution.

The results from the third- and fourth-order SF-PIF methods
produce well-acceptable results compared to
the corresponding RK methods. 
Except that there are minor differences
in the shape of Kelvin-Helmholtz instabilities
along the primary slip line and the bottom jet,
the overall dynamics of the two SF-PIF solutions 
match well with the RK solutions, validating
the fidelity of the proposed SF-PIF methods
in the presence of a strong shock.


\subsection{3D Euler equations}\label{sec:results-3d}

\subsubsection{3D Sedov test}\label{sec:3dsedov}

To test the code's ability to maintain the spherical
symmetry in all spatial directions, 
we consider the Sedov blast test~\cite{sedov1993similarity} in 3D.
Initially, a point-source of a highly pressurized perturbation is given
at the domain center, which leads to a strong spherical shock wave propagating
outward from the source. The simulation domain is a 3D square box of
\( [-1.2, 1.2] \times [-1.2, 1.2] \times [-1.2, 1.2] \) 
resolved on a \( 128 \times 128 \times 128 \) grid resolution.
Outflow boundary conditions are imposed in all directions.
The initial conditions are based on the setup found in~\cite{fryxell2000flash},
but we choose the deposited energy, \( E_{\text{tot}} = 0.851072 \),
according to the setup in~\cite{boscheri2014direct}.

\begin{figure}[hb!]
    \centering
    \begin{subfigure}{80mm}
        \centering
        \includegraphics[width=0.95\textwidth]{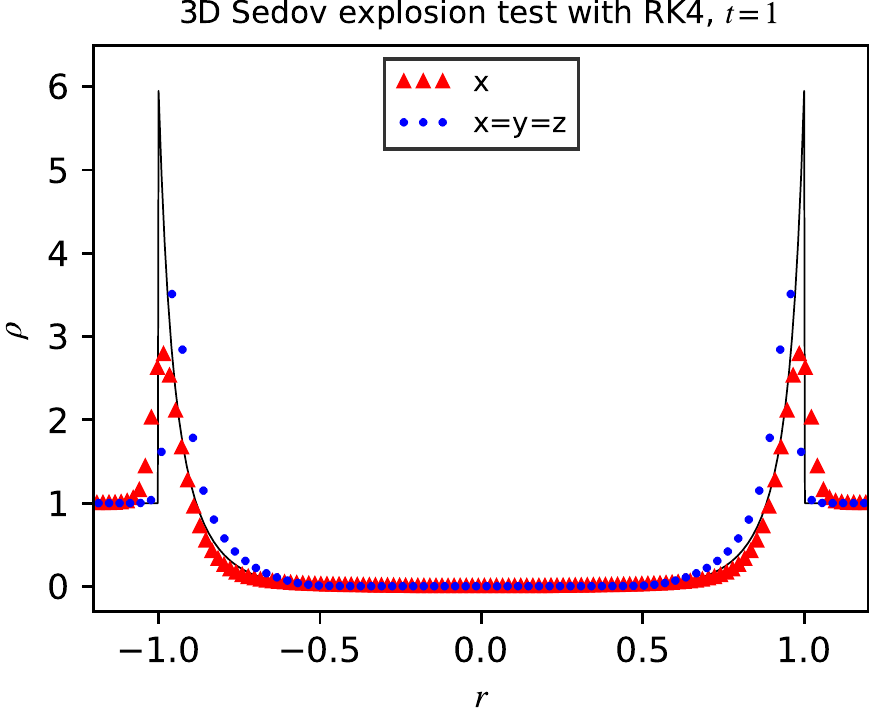}
    \end{subfigure}
    \begin{subfigure}{80mm}
        \centering
        \includegraphics[width=0.95\textwidth]{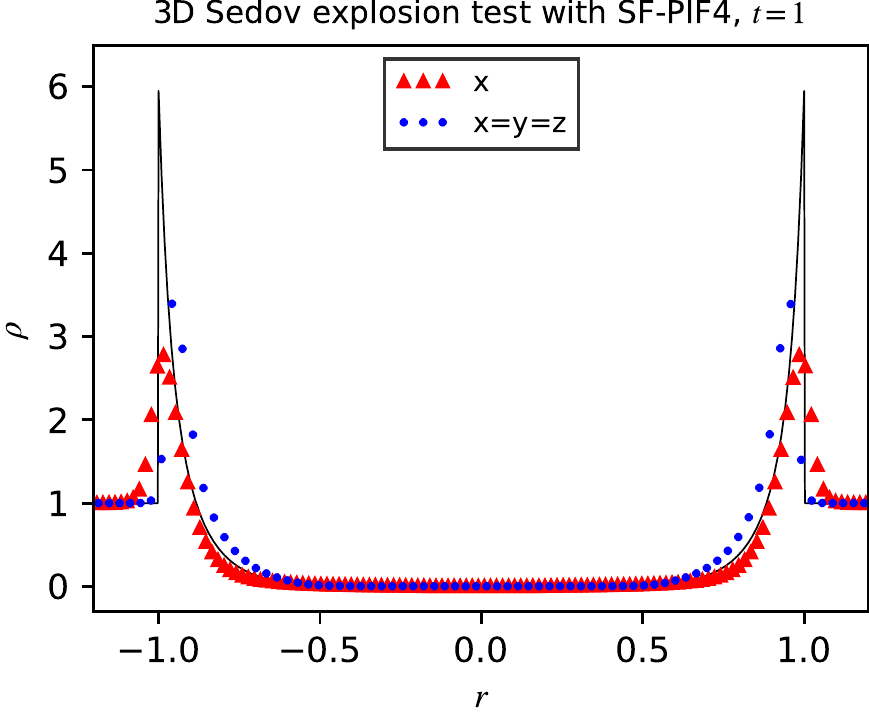}
    \end{subfigure}
    \caption{The density profiles of the 3D Sedov explosion test at \( t = 1 \).
        The red triangles represent the density profiles along the
        \( x \)-axis from the origin, while the blue circles represent
        the density profiles along the diagonal axis.
        The solid black line is the well-known exact self-similar solution.
        All simulations are performed on a \( 128 \times 128 \times 128 \) grid resolution,
        solved with RK4 (\textbf{left}) and SF-PIF4 (\textbf{right}) temporal solvers.
    }\label{fig:3dsedov}
\end{figure}

\cref{fig:3dsedov} shows the density profiles at \( t = 1 \) along the diagonal \( (x = y = z) \)
and the \( x \)-axis \( (y = z = 0) \). The results on the left panel are solved with the
RK4 method and the right panel with SF-PIF4. The exact self-similar solution
is plotted as a reference solution in solid black curves on both figures.
As shown, both the RK4 and SF-PIF4 methods show
nearly identical results. Despite the visible differences between
the diagonal and the \( x \)-axis profiles, especially at the highest
peak and the shock front location, both fourth-order temporal solvers
show well-maintained reflectional symmetry along the normal axis at the origin.

\subsubsection{3D Sod's blast explosion test}\label{sec:3dsod}
\begin{figure}
    \centering
    \begin{subfigure}{80mm}
        \centering
        \includegraphics[width=0.95\textwidth]{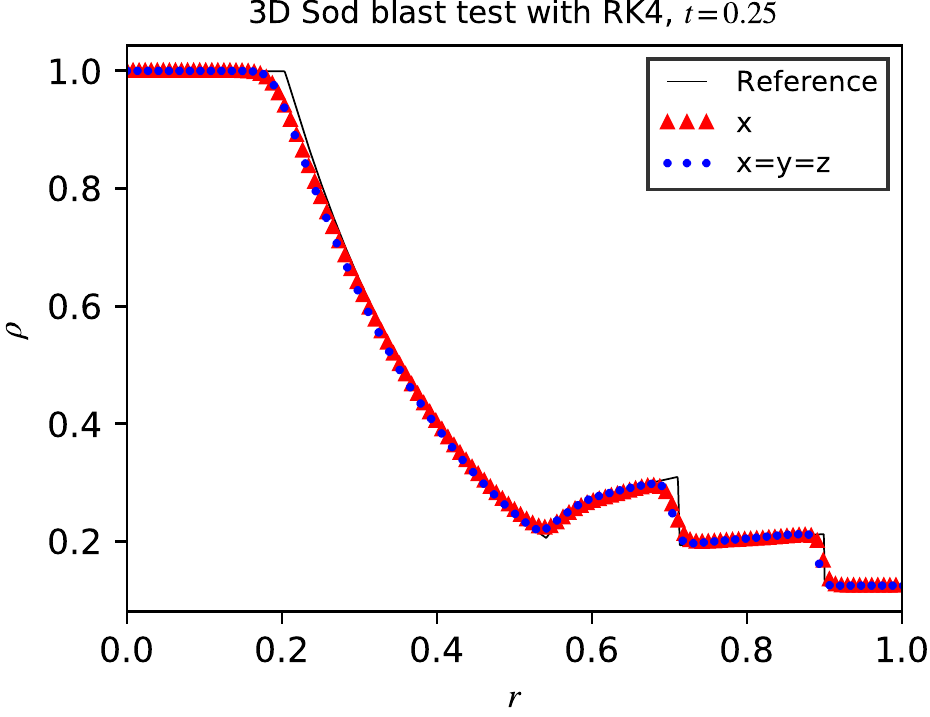}
    \end{subfigure}
    \begin{subfigure}{80mm}
        \centering
        \includegraphics[width=0.95\textwidth]{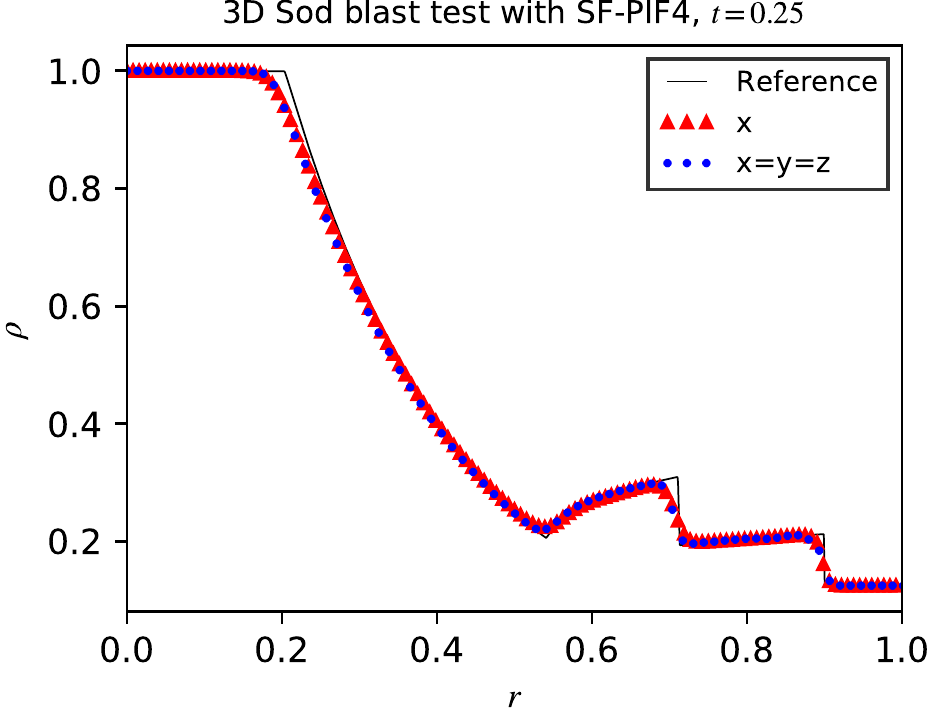}
    \end{subfigure}
    \caption{The density profiles of the 3D Sod's blast test at \( t = 0.25 \).
        The red triangles represent the density profiles along the
        \( x \)-axis from the origin, and the blue circles represent
        the density profiles along the diagonal axis.
        The solid black curve is the reference solution
        calculated with the 1D Euler equations with the geometric source term.
        Simulations are performed on a \( 256 \times 256 \times 256 \) grid resolution,
        solved with RK4 (\textbf{left}) and SF-PIF4 (\textbf{right}) temporal solvers.
    }\label{fig:3dsod}
\end{figure}

Next, we consider the 3D explosion test problem found in~\cite{boscheri2014direct}.
This test is a three-dimensional extension of 
the 1D Sod's shock tube problem~\cite{sod1978survey}, 
which we already solved in a rotated 2D shock-tube configuration in~\cref{sec:sod45}.
The calculations are performed on a \( [-1, 1] \times [-1, 1] \times [-1, 1] \)
domain with outflow boundary conditions using a \( 256 \times 256 \times 256 \)
grid resolution.

The result of the density profiles along the diagonal \( (x = y = z) \)
and \( x \)-axis \( (y = z = 0) \) at \( t = 0.25 \) are presented in~\cref{fig:3dsod}.
The solid black curve represents the reference solution
using the 1D Euler equations with the appropriate geometric source term
according to~\cite{boscheri2014direct}.
The results show that the SF-PIF4 solver captures the shock
profile and the spherical symmetry exceedingly well compared to
the RK4 result and the reference 1D solution.
There is no noticeable difference between
the two fourth-order temporal solvers.

\subsubsection{3D Riemann problem}\label{sec:3drp}

Finally, we performed the 3D Riemann problem presented in~\cite{balsara2015three}.
We follow the same initial conditions, consisting of 
eight constant initial conditions
in each octant of the computational domain,
\( [-1, 1] \times [-1, 1] \times [-1, 1]\),
resolved on a \( 256 \times 256 \times 256 \) grid resolution.
The setup imposes outflow boundary conditions at all boundaries.
The initial condition will carry out 2D Riemann problems at each octant interface,
including the diagonal plane of the 3D computational cubic.

\begin{figure}
    \centering
    \begin{subfigure}{160mm}
        \centering
        \includegraphics[width=0.95\textwidth]{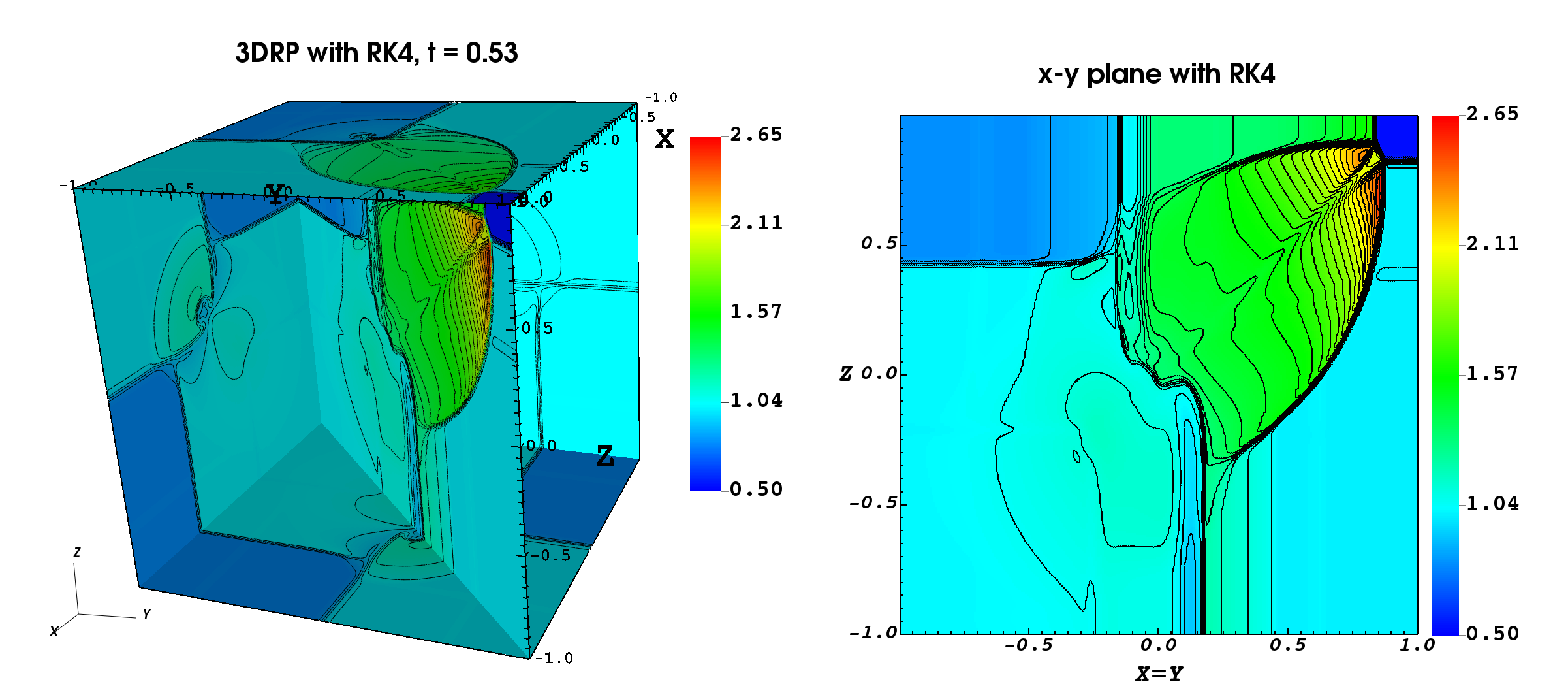}
    \end{subfigure}
    \begin{subfigure}{160mm}
        \centering
        \includegraphics[width=0.95\textwidth]{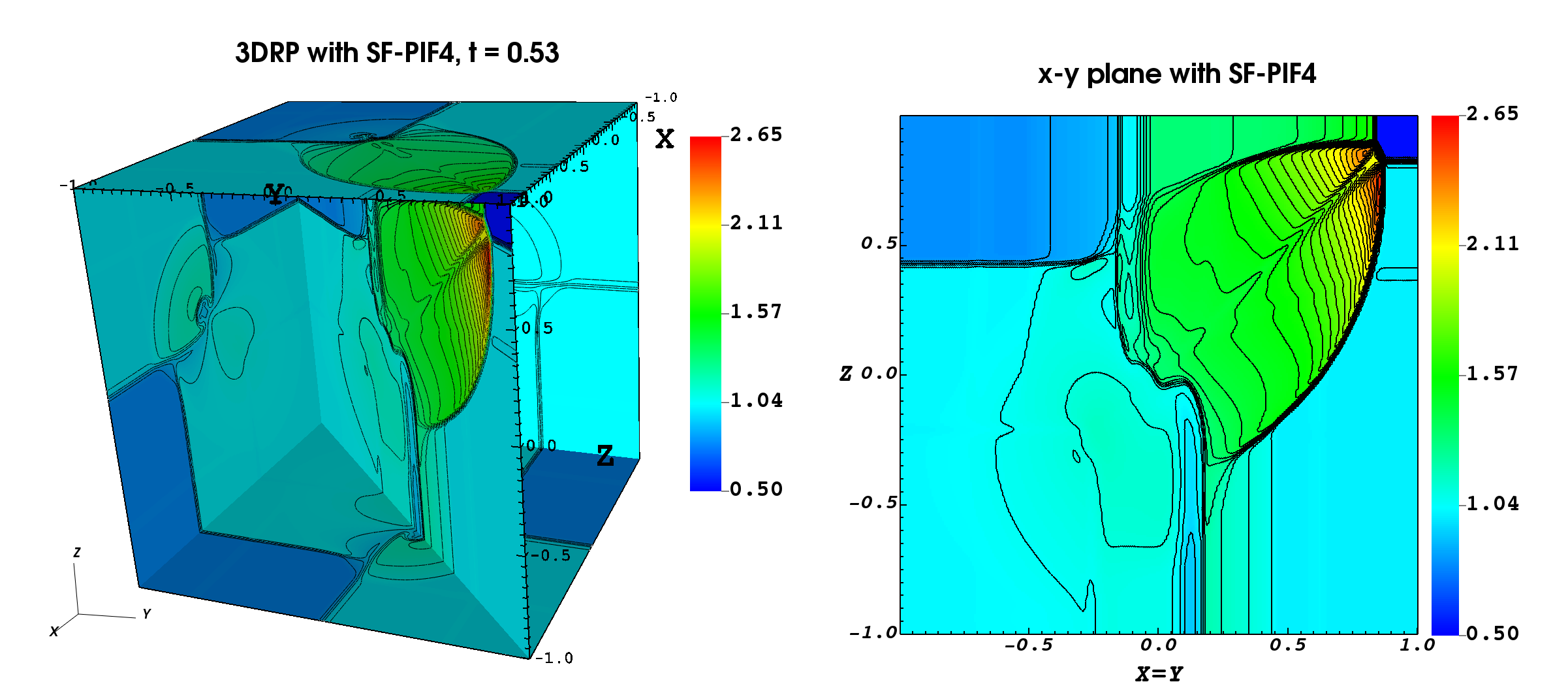}
    \end{subfigure}
    \caption{The density maps of the 3D Riemann problem test at \( t = 0.53 \).
        Forty contour lines are over-plotted.
        The left panels show each face's geometrical views, while 
        the right panels show
        the detailed picture of the diagonal planes.
        All simulations are performed on a \( 256 \times 256 \times 256 \) grid resolution,
        solved with RK4 (\textbf{top}) and SF-PIF4 (\textbf{bottom}) solvers.
    }\label{fig:3drp}
\end{figure}

The resulting density profiles at \( t = 0.53 \) are given in~\cref{fig:3drp}.
The pseudo-color map ranges between \( [0.5, 2.65] \),
and we over-plot 40 levels of contour lines using the same range.
We observe three different 2D Riemann problems on the
left, top, and the diagonal planes in the left panel. The diagonal
planes are separately shown in the right panel.
SF-PIF4 method is able to capture all the important features
as much as the RK4 result, confirming the validity of the SF-PIF4 method
in comparison.

\begin{table}[ht!]
    \footnotesize
    \centering
    \caption{Performance results for the 3DRP test problem.
        All performance results (measured in seconds) are averaged over
        five simulation runs conducted on 16 nodes of \textit{lux} cluster.
        Each node has 2 \( \times \) 20-core Intel Xeon Gold 6248 (Cascade Lake) CPUs,
        and we utilized 512 parallel threads for each run.
    }\label{table:perform-3drp}
    \begin{tabular}{@{}cccccccccccc@{}}
        \toprule
        \multirow{2}{*}{Grid Resolution} & \multicolumn{2}{l}{RK3} &  & \multicolumn{2}{l}{SF-PIF3}
                & & \multicolumn{2}{l}{RK4} & & \multicolumn{2}{l}{SF-PIF4} \\
        \cmidrule(lr){2-3} \cmidrule(l){5-6} \cmidrule(l){8-9} \cmidrule(l){11-12}
        & CPU Time & Speedup &  & CPU Time & Speedup & & CPU Time & Speedup & & CPU Time & Speedup\\
        \midrule
        \( 64 \times 64 \)   & \SI{1.79}{\second} & 1.0 &  & \SI{1.09}{\second} & 0.61 & &
            \SI{2.95}{\second} & 1.64 &  & \SI{2.67}{\second} & 1.49 \\
        \( 128 \times 128 \) & \SI{15.62}{\second} & 1.0 &  & \SI{7.63}{\second} & 0.49 & &
            \SI{25.88}{\second} & 1.66 &  & \SI{18.97}{\second} & 1.21 \\
        \( 256 \times 256 \) & \SI{191.00}{\second}   & 1.0 &  & \SI{82.02}{\second} & 0.43 & &
            \SI{321.34}{\second} & 1.68 &  & \SI{201.40}{\second} & 1.05 \\
        \( 512 \times 512 \) & \SI{2679.85}{\second}  & 1.0 &  & \SI{1173.54}{\second} & 0.44 & &
            \SI{4507.88}{\second} & 1.68 &  & \SI{2817.78}{\second} & 1.05 \\
    \end{tabular}
\end{table}

\cref{table:perform-3drp} shows the performance results for the 3D Riemann problem test
on four grid resolutions. As shown in the table, the SF-PIF4 method demonstrates nearly
the same performance as the \textit{third-order} RK method,
especially in high-resolution cases.
We should note that the performance gains from the SF-PIF methods are
more compensated on the high-resolution grids,
which are indispensable for high fidelity physical simulation studies.

\section{Conclusion}\label{sec:conclusion}


In this study, we have extended the original third-order (non-recursive) SF-PIF method~\cite{lee2021single} to a fourth-order temporal scheme 
with the improved recursive version of the system-free (SF) approach. The newly proposed recursive SF approach enables the fourth-order extension of the SF-PIF method (SF-PIF4),
increasing computational performance gain with a simplistic code structure.

The original SF approach's critical design purpose~\cite{lee2021single} is to bypass all the analytical derivations of Jacobians, Hessians, and even higher-order derivatives tensor terms. With the SF approach, approximating the tensor contractions of \textit{Jacobian-like} terms in the 
Lax-Wendroff type time discretization becomes simplified.

In this paper, the SF method's advantage is further enhanced by introducing a new recursive procedure.  The recursive SF method requires only a small amount of calculations for approximating the tensor contractions, thereby empowering the fourth-order extension of the SF-PIF method to ease and faster performance.

We have tested our fourth-order and third-order SF-PIF methods with a wide range of test problems in two and three spatial dimensions. The results show the SF-PIF methods have significantly faster computational time performance  than the corresponding SSP-RK methods at the same temporal order while maintaining the equivalent solution accuracy. In two-dimensional cases, SF-PIF methods show more than two times faster performance results than the SSP-RK counterparts.  Moreover, we demonstrate that the fourth-order SF-PIF's performance results are nearly the same as the \textit{third-order} SSP-RK method.
This enhanced performance gain in our SF-PIF solvers can provide a big leap in large-scale parallel computing to save computational costs and reach highly accurate numerical predictions in practice. 

We believe the recursive SF method has a broad potential to be relevant in various fields of study. It is neither designed particularly for the PIF method nor any specific numerical methods in computational fluid dynamics. Instead, it is solely intended for approximating the \textit{Jacobian-like} tensor contractions. We presume that our recursive SF schemes are applicable in other numerical algorithms to enhance the calculation speed and ease the code implementations wherever the Jacobians, Hessians, or higher derivative tensors are required.

A further extension is to design the SF-PIF method in arbitrary order. As reported in~\cref{sec:vortex}, the temporal errors dominate the spatial errors in high-resolution simulations; thus, it is noteworthy to use the same accuracy in both the spatial and temporal solvers. We realize that a naive extension of the SF-PIF method to the fifth or higher-order could potentially be less attractive due to the drastic increase of complexity in Taylor expansion. We aim to reduce such complexity and make the SF-PIF method an arbitrary order of accuracy, which will be further investigated in our future studies.

\section{Acknowledgements}
We acknowledge the use of the lux supercomputer at
UC Santa Cruz, funded by NSF MRI grant AST 1828315.

\bibliography{refs}\label{sec:references}

\end{document}

%% file: header.tex
\usepackage{framed,multirow}

\usepackage{amssymb}
\usepackage{latexsym}

\usepackage{lineno}
\usepackage[bookmarks=true,bookmarksnumbered=false,bookmarksopen=false,breaklinks=false,pdfborder={0 0 0},pdfborderstyle={},backref=false,colorlinks=true]{hyperref}
\usepackage{color}
\usepackage{amsfonts}
\usepackage{amsmath}
\usepackage{amssymb}
\usepackage{multirow}
\usepackage{graphicx}
\usepackage{rotating}
\usepackage{mathtools}
\usepackage[shortlabels]{enumitem}

\usepackage{siunitx} 
\usepackage{graphicx}
\usepackage{subcaption}
\usepackage[capitalize]{cleveref}

\usepackage{booktabs}
\usepackage{multirow}

\newcommand{\bff}{\mathbf{f}}

\newcommand{\half}{\frac{1}{2}}

\newcommand{\mm}[1]{\rm mm}

\newcommand{\beq}{\begin{equation}}
\newcommand{\eeq}{\end{equation}}
\newcommand{\bea}{\begin{eqnarray}}
\newcommand{\eea}{\end{eqnarray}}

\newcommand{\dt}{\Delta t}
\newcommand{\dx}{\Delta x}
\newcommand{\dy}{\Delta y}
\newcommand{\dz}{\Delta z}

\newcommand{\bit}{\begin{itemize}}
\newcommand{\eit}{\end{itemize}}
\newcommand{\ben}{\begin{enumerate}}
\newcommand{\een}{\end{enumerate}}

\newcommand{\bF}{\mathbf{F}}
\newcommand{\bG}{\mathbf{G}}
\newcommand{\bH}{\mathbf{H}}

\newcommand{\bU}{\mathbf{U}}
\newcommand{\bV}{\mathbf{V}}
\newcommand{\bW}{\mathbf{W}}

\newcommand{\bX}{\mathbf{X}}

\newcommand{\bg}{\mathbf{g}}
\newcommand{\bh}{\mathbf{h}}

\newcommand{\bx}{\mathbf{x}}

\newcommand{\Div}{\nabla^{f}}
\newcommand{\scolon}{\,;\,}

\let\oldequation\equation
\let\oldendequation\endequation
\renewenvironment{equation}
  {\linenomathNonumbers\oldequation}
  {\oldendequation\endlinenomath}

\makeatletter
\def\@xfootnote[#1]{%
  \protected@xdef\@thefnmark{#1}%
  \@footnotemark\@footnotetext}
\makeatother

\usepackage{url}
\usepackage{xcolor}
\definecolor{newcolor}{rgb}{.8,.349,.1}